\documentclass[11pt, reqno]{amsart}
\usepackage{amssymb}
\usepackage{amsfonts}
\usepackage{amsmath}
\usepackage{mathtools}
\usepackage{amsthm}
\usepackage{xcolor}
\usepackage{hyperref}
\usepackage{amsaddr}
\usepackage{enumitem}
\usepackage{graphicx}
\usepackage{bmpsize}
\usepackage{ifthen}
\newboolean{images}
\newboolean{lowres}
\setboolean{images}{true}
\setboolean{lowres}{true}
\usepackage{comment}

\newtheorem{theorem}{Theorem}[section]
\newtheorem{proposition}[theorem]{Proposition}
\newtheorem{definition}[theorem]{Definition}
\newtheorem{lemma}[theorem]{Lemma}
\newtheorem{corollary}[theorem]{Corollary}
\newtheorem{example}[theorem]{Example}
\newtheorem{conjecture}[theorem]{Conjecture}
\newtheorem{remark}[theorem]{Remark}
\newtheorem{problem}[theorem]{Problem}

\usepackage{geometry}
\def\C{\mathcal{C}}

\def\Z{\mathbb{Z}}

\def\P{\mathcal{P}}

\DeclareMathOperator{\dev}{dev}
\DeclareMathOperator{\Aut}{Aut}
\DeclareMathOperator{\Atop}{Atop}

\DeclareMathOperator{\Mult}{Mult}
\DeclareMathOperator{\supp}{supp}

\newcommand{\olsi}[1]{\,\overline{\!{#1}}} % overline short italic

\begin{document}

\title{Three-dimensional symmetric designs of propriety~3}

\author[A.~Bahmanian, V.~Kr\v{c}adinac, L.~Reli\'{c}, and S.~Suda]{Amin Bahmanian$^1$, Vedran Kr\v{c}adinac$^2$, Lucija Reli\'{c}$^2$, and Sho Suda$^3$}

\address{$^1$Department of Mathematics, Illinois State University,  USA}

\address{$^2$University of Zagreb Faculty of Science, Department of Mathematics, Bijeni\v{c}ka cesta~30, 10000~Zagreb, Croatia}

\address{$^3$Department of Mathematics, National Defense Academy of Japan, 1-10-20 Hashirimizu, Yokosuka, Kanagawa 239-8686, Japan}

\email{mbahman@ilstu.edu}
\email{vedran.krcadinac@math.hr}
\email{lucija.relic@math.hr}
\email{ssuda@nda.ac.jp}



\date{September 21, 2026}

\begin{abstract}
We define symmetric designs of dimension~$n$ and propriety~$d$, providing a unifying generalization of several classes of higher-dimensional symmetric designs previously studied. We focus on the case $n=d=3$, which leads to the following question: Can we fill the $v^3$ cells of a $v\times v\times v$ cube with $\{0,1\}$ in such a way that each layer parallel to each face contains a fixed number~$k$ of ones, and that for every two parallel layers there are exactly $\lambda$ positions where they have matching ones? We establish necessary conditions on the parameters $(v,k,\lambda)$,  introduce notions of difference sets and multipliers for these objects, and enumerate small examples up to equivalence. Furthermore, we construct infinite families of these objects using difference sets, symmetric designs, doubly regular tournaments, Hadamard matrices, Latin cubes,
and association schemes on triples.
\end{abstract}

\maketitle

\allowdisplaybreaks

\section{Introduction}

In light of recent interest in higher-dimensional analogues of various combinatorial objects such as permutations \cite{LL14,LL16}, Latin squares \cite{BSid23, BProc24, BArxiv22, MW08}, Hadamard matrices \cite{BS, KPT23}, association schemes \cite{BB23, BB25, BP21}, and symmetric block designs \cite{KP25, KPT25, KR25}, we consider the following problem. Throughout the paper, we assume $v\ge 2$ and $2\le d\le n$.

\begin{problem}\label{prob1}
    Given integers $v\ge 2$ and $k,\lambda\ge 0$, can we fill the $v^3$ cells of a $v\times v\times v$ cube with $\{0,1\}$ in such a way that each layer parallel to each face contains  $k$ ones, and that for every two parallel layers there are exactly $\lambda$ positions where they have matching ones?
\end{problem}

To elaborate, let us provide some definitions first.
A \emph{Hadamard matrix} of order~$v$ is a $v\times v$ matrix over $\{1,-1\}$ such
that its rows are mutually orthogonal, as well as its columns. Shlichta~\cite{PS79}
defined an \emph{$n$-dimensional Hada\-mard matrix} of order~$v$ and propriety~$d$ as
a $v\times \cdots \times v$ array~$H$ of dimension~$n$ over $\{1,-1\}$ such that
all $(d-1)$-dimensional parallel subarrays of~$H$ are orthogonal. The definition
can also be found in~\cite[Chapter~5]{KH07} and~\cite[Chapter~6]{YNX10}.

A \emph{Latin square} of order~$v$ is a $v\times v$ matrix over $\{1,\ldots,v\}$
such that each symbol appears exactly once in every row and column. A
\emph{Latin hypercube} of dimension~$n$, order~$v$, and propriety~$d$ can be defined
as a $v\times \cdots \times v$ array~$L$ of dimension~$n$ over $\{1,\ldots,v\}$
such that each symbol appears equally often in every $(d-1)$-dimensional
subarray of~$L$, namely $v^{d-2}$ times. This definition can be found
in~\cite[Chapter~3]{LM98}, where \emph{type} $n-d+1$ is used instead
of propriety~$d$.

In both cases, objects of propriety~$d<n$ are also of higher propriety
$d+1,\ldots,n$. The strongest possible propriety is $d=2$, meaning that
all $2$-dimensional subarrays of~$H$ are Hadamard matrices, and of~$L$
Latin squares. Objects of propriety~$2$ are called \emph{proper}
$n$-dimensional Hadamard matrices, or $n$-dimensional
\emph{permutation cubes}~\cite{KD15}.
De~Launey~\cite{WdL90} generalized the propriety~$2$ definition to a
wide variety of combinatorial designs, including symmetric block designs,
orthogonal designs, generalized Hadamard matrices, and generalized weighing
matrices. This generalization is called \emph{proper $n$-di\-men\-sio\-nal designs}
in~\cite[Definition~2.7]{WdL90}.

An incidence matrix of a symmetric $(v,k,\lambda)$ design is a $v\times v$ matrix
over $\{0,1\}$ such that every row and column contains exactly~$k$ entries~$1$, and scalar
products of all pairs of distinct rows and columns are~$\lambda$. We remark that constructing an infinite family of symmetric designs or even a single symmetric design is a notoriously difficult problem and, for example, the existence of finite projective planes of non-prime power orders is a long-standing open problem~\cite{BJL99, IS06, IT07, ESL83}. Symmetric designs of dimension~$n$ and
propriety~$2$ have recently been studied in~\cite{KP25, KPT25} under the name
\emph{cubes of symmetric designs} or \emph{$\C^n(v,k,\lambda)$-cubes}. These are
$v\times\cdots\times v$ arrays of dimension~$n$ whose $2$-dimensional
subarrays are incidence matrices of $(v,k,\lambda)$ designs.

We now generalize $n$-dimensional symmetric designs to arbitrary
propriety~$d$. To make the definition and the previous discussion clear,
let us note that an \emph{$n$-dimensional matrix} or \emph{array}~$A$ of
order~$v$ over a set of symbols~$S$ is a function $A:\{1,\ldots,v\}^n \to S$.
The domain is the Cartesian product of $n$ copies of $\{1,\ldots,v\}$.
A $d$-dimensional \emph{subarray} is obtained by fixing $n-d$ coordinates.
Two such subarrays are considered \emph{parallel} if the same coordinates
are fixed and the values of the fixed coordinates agree, except one.

\begin{definition}\label{maindef}
A \emph{symmetric $(v,k,\lambda)$ block design of dimension~$n$ and
propriety~$d$} is a $v\times\cdots\times v$ array~$A$ of dimension~$n$ over
$\{0,1\}$ such that every $(d-1)$-dimensional subarray of~$A$ contains
exactly~$k$ entries~$1$, and scalar products of all pairs of parallel
$(d-1)$-dimensional subarrays are~$\lambda$. The set of all such objects
will be denoted by $\C^n_d(v,k,\lambda)$.
\end{definition}

Some direct consequences of the definition follow. First, cubes of symmetric
designs from~\cite{KP25, KPT25} are the propriety~$2$ case:
$$\C^n(v,k,\lambda)= \C^n_2(v,k,\lambda).$$
Next, $n$-dimensional symmetric designs of propriety~$d<n$ are also of
higher propriety, but with different parameters:
$$\C^n_d(v,k,\lambda) \subseteq \C^n_{d+1}(v,vk,v\lambda).$$
Indeed, every $d$-dimensional subarray is the disjoint union of $v$ parallel $(d-1)$-dimensional subarrays, and the scalar product of two parallel $d$-dimensional subarrays is the sum of the $v$ corresponding scalar products.
Consequently, $\C^n(v,k,\lambda)\subseteq \C^n_d(v,v^{d-2}k,v^{d-2}\lambda)$
for all $d=2,\ldots,n$; for $d>2$, this inclusion can be strict. In
Section~\ref{sec3}, we give examples of $\C^3_3(v,vk,v\lambda)$-cubes
that are not in $\C^3(v,k,\lambda)$.

Another generalization of symmetric designs to higher dimension was
introduced in~\cite{KR25} and studied in~\cite{KP25} under the name
\emph{projection cubes} or \emph{$\P^n(v,k,\lambda)$-cubes}. For an
$n$-dimensional array~$A$, a projection~$\Pi_{ij}(A)$ is the
$2$-dimensional matrix with $(x_i,x_j)$-entry
\begin{equation}\label{projsum}
\sum_{1\le x_1,\ldots,\widehat{x_{i}},\ldots,\widehat{x_{j}},\ldots,x_n\le v} A(x_1,\ldots,x_n).
\end{equation}
The sum is taken over all $n$-tuples $(x_1,\ldots,x_n)\in \{1,\ldots,v\}^n$
with fixed coordinates $x_i$ and $x_j$ in a field of characteristic~$0$.
The defining property of $\P^n(v,k,\lambda)$-cubes is that all projections
$\Pi_{ij}(A)$, $1\le i<j\le n$ are incidence matrices of symmetric $(v,k,\lambda)$
designs. Projection cubes are reminiscent of Room cubes~\cite[Section~VI.50.8]{JHD07}
and Costas cubes~\cite{JY18, USO25}. For $n\ge 3$, they are also a special case of Definition~\ref{maindef}:
\begin{equation}\label{projincl}
\P^n(v,k,\lambda) \subseteq \C^n_n(v,k,0).
\end{equation}
Fixing $x_i$ and summing~\eqref{projsum} over $x_j\in\{1,\ldots,v\}$ gives
the row sum of a $(v,k,\lambda)$ incidence matrix, which equals~$k$.
This shows that every $(n-1)$-dimensional subarray of an array
$A\in\P^n(v,k,\lambda)$ contains exactly~$k$ entries~$1$. The condition
that parallel $(n-1)$-dimensional subarrays are orthogonal
follows because the sums~\eqref{projsum} are in $\{0,1\}$.
A pair of non-orthogonal parallel subarrays would have a common
$1$-entry, and we could find a sum~\eqref{projsum} that is at
least~$2$. Again, in Section~\ref{sec3} we give examples
of $\C^3_3(v,k,0)$-cubes that are not in $\P^3(v,k,\lambda)$,
showing that the inclusion~\eqref{projincl} is strict.

Our final observation for arbitrary dimension~$n$ and propriety~$d$ concerns complementation. If $A$ is a $\C^n_d(v,k,\lambda)$-cube, then $J-A$ is a $\C^n_d(v,v^{d-1}-k,v^{d-1}-2k+\lambda)$-cube. Here and in the sequel, $J$ denotes the $n$-dimensional array with all entries~$1$, so that $J-A$ is the array~$A$ with reversed $0$- and $1$-entries.

In the rest of the paper we study the special case $n=d=3$, i.e.\
three-dimensional symmetric designs of propriety~$3$. Problem~\ref{prob1} is
the existence question for this case. As already mentioned,
every $\C^3(v,k',\lambda')$-cube from~\cite{KP25, KPT25} is an example
of a $\C^3_3(v,k,\lambda)$-cube for $k=vk'$ and $\lambda=v\lambda'$, but
this is a very restricted subfamily. For one, the parameters~$k$
and~$\lambda$ need not be divisible by~$v$. Moreover, the
equality $\lambda'(v-1)=k'(k'-1)$ is a necessary condition only
for propriety $d=2$. For cubes of propriety $d=3$, the
parameter~$\lambda$ need not be uniquely determined by~$v$ and~$k$.
Examples for both situations appear in Tables~\ref{tab1}--\ref{tab4}.

The paper is organized as follows. In Section~\ref{sec2}, we prove necessary
conditions for the existence of $\C^3_3(v,k,\lambda)$-cubes in the form of
inequalities. The extremal cases are studied and characterized. For small
parameters~$v$ and~$k$, all $\C^3_3(v,k,\lambda)$-cubes are enumerated up
to equivalence.

Section~\ref{sec3} contains results about difference sets
for $\C^3_3(v,k,\lambda)$-cubes. Cubes constructed from
difference sets are characterized as possessing an autotopy
group acting sharply transitively on the three coordinates.
Small examples of difference sets are also enumerated.
Multipliers are defined and used to push the enumeration
a little further. Constructions starting from ordinary
difference sets of Paley and Hadamard types are proved,
providing infinitely many examples of difference sets
for $\C^3_3(v,k,\lambda)$-cubes.

A noteworthy aspect of this paper is the way various combinatorial
objects are interrelated and how their connections give rise to
infinite families of three-dimensional symmetric designs.
In Section~\ref{sec4}, we provide constructions for $\C^3_3(v,k,\lambda)$-cubes
summarized in the following theorem. We  show by examples that the constructions
give rise to multiple nonequivalent cubes with the same parameters.

\begin{theorem} \label{infifamthmsymcub}
   \begin{enumerate}[label=\textup{(\roman*)}]%[label=\textup{\arabic*.}]
   \item []
       \item If there is a symmetric $(v,k,\lambda)$ design, then there exist
       $\C_3^3(v,vk,v\lambda)$-cubes, and
       $\C_3^3(v,K_i,\Lambda_i)$-cubes, $i\in\{0,1\}$, where
\begin{alignat*}{3}
K_0 &=k (3 v + 4\lambda -6 k ), & \quad \Lambda_0 &=k ( v - 4 \lambda-2 ) + 2 \lambda (2 \lambda+1),\\[1mm]
K_1 &=k (3 v + 4\lambda -2 k) - 4 v\lambda, & \quad \Lambda_1 &= k (v - 4 \lambda+2) + 2 \lambda (2 \lambda-1).
\end{alignat*}
\item If there is a doubly regular tournament of order $4t-1$, then there exist $\C_3^3(4t-1,K_i,\Lambda_i)$-cubes, $i\in\{0,1\}$, where
\begin{alignat*}{3}
K_0 &=(2t-1)(4t+3), & \quad \Lambda_0 &=(2t-1)(2t+3),\\[1mm]
K_1 &=2(4t^2-5t+2), & \quad \Lambda_1 &= 4(t-1)^2.
\end{alignat*}
        \item If there is a Hadamard matrix of order $v^2$, then there exist cubes in
\[ \C^3_3(v^3,\, v^4(v^2-1)/2,\, v^4(v^2-2)/4). \]
       \item  If there is  a regular Hadamard matrix of order $v^2$, then there exist cubes in
\[ \C^3_3(v^3,\, v^3(v^3-1)/2,\, v^3(v^3-2)/4). \]
       \item  If there is a symmetric association scheme on triples of order~$v$, then there exists a family of $\C^3_3(v,k,\lambda )$-cubes, where the parameters
       $k$ and $\lambda$ are determined by the intersection numbers of the association
       scheme on triples.
   \end{enumerate}
\end{theorem}
Constructions~(iii) and~(iv) are used to obtain three-dimensional symmetric designs admitting further regularities. Additionally,
%inspired by Bruck's lemma~\cite{RHB55} and its generalization by Kantor~\cite{Ka69},
we embed three-dimensional symmetric designs in other three-dimensional symmetric designs.
Finally, in Section~\ref{sec5}, we give some concluding remarks and indicate directions for further research.

\section{Necessary conditions and small examples}\label{sec2}

For $n=d=3$, the conditions of Definition~\ref{maindef} can be written as
\begin{align}\label{maincond}
\begin{cases}
\displaystyle\sum_{y=1}^v \sum_{z=1}^v A(x_1,y,z)A(x_2,y,z) = [x_1=x_2]\, k + [x_1\neq x_2]\,\lambda, &\mbox{ for } x_1,x_2\in \{1,\dots,v\}, \\
\displaystyle\sum_{x=1}^v \sum_{z=1}^v A(x,y_1,z)A(x,y_2,z) = [y_1=y_2]\, k + [y_1\neq y_2]\,\lambda, &\mbox{ for } y_1,y_2\in \{1,\dots,v\}, \\
\displaystyle\sum_{x=1}^v \sum_{y=1}^v A(x,y,z_1)A(x,y,z_2) = [z_1=z_2]\, k + [z_1\neq z_2]\,\lambda,  &\mbox{ for } z_1,z_2\in \{1,\dots,v\}.
\end{cases}
\end{align}
The square brackets are the Iverson symbol: $[P]=1$ if the statement~$P$
is true, and $[P]=0$ otherwise (see~\cite{DK92}).
Obvious necessary conditions for the existence of
$\C^3_3(v,k,\lambda)$-cubes are the inequalities
$0\le k \le v^2$ and $0\le \lambda \le k$.
We first prove an inequality involving parameters of cubes
of propriety~$3$, analogous to the necessary condition
$\lambda'(v-1)=k'(k'-1)$ for cubes of propriety~$2$. Throughout the
paper, we shall refer to $1$-dimensional and $2$-dimensional subarrays
as \emph{rows} and \emph{layers}, respectively.

\begin{proposition}\label{prop:1}
If a cube $A\in \C^3_3(v,k,\lambda)$ exists, then
\begin{equation}\label{leqnec}
\lambda v(v-1)\geq k(k-v).
\end{equation}
Equality holds if and only if every row-sum of~$A$ is $k/v$.
In that case, $(v,k,\lambda)=(v,vk',v\lambda')$ holds,
where $(v,k',\lambda')$ are admissible parameters of cubes of
propriety~$2$.
\end{proposition}

\begin{proof}
Let $s_{xy}=\sum_{z=1}^v A(x,y,z)$ be a row-sum of~$A$ in the direction of the $z$-axis,
for $x,y\in \{1,\dots, v\}$. By~\eqref{maincond}, we have $\sum_{x=1}^v s_{xy} = \sum_{y=1}^v s_{xy} = k$, and
\begin{equation}\label{sumentries}
\sum_{x=1}^v \sum_{y=1}^v s_{xy} = vk.
\end{equation}
We count in two ways the size of the set
$$S:=\Big\{\big(x,y,\{z_1,z_2\}\big) \kern 1mm\mathop{\big|}\kern 1mm A(x,y,z_1)=A(x,y,z_2)=1,\kern 1mm x,y,z_1,z_2\in \{1,\ldots,v\}, \kern 1mm z_1<z_2 \Big\}.$$
On the one hand, since parallel $z$-layers share $\lambda$ common
$1$-entries, $|S|= \lambda \binom{v}{2}$, and on the other hand,
$|S|=\sum_{x=1}^v \sum_{y=1}^v \binom{s_{xy}}{2}$. The latter sum reaches a minimal value under the constraint~\eqref{sumentries}
when the variables are equal. Hence, $\binom{v}{2}\lambda \ge v^2 \binom{k/v}{2}$ and
equality holds if and only if $s_{xy}=k/v$ for all $x$ and $y$. Since the argument is
also valid for row-sums in the direction of the $x$-axis and $y$-axis, we conclude that
equality holds if and only if every row-sum of~$A$ is $k/v$. Denoting the common row-sum by~$k'$
and substituting $k=vk'$ in the equality gives $\lambda=\frac{vk'(k'-1)}{v-1}$.
Since $v$ and $v-1$ are relatively prime, $v-1$ must divide $k'(k'-1)$,
i.e.\ $\lambda=v\lambda'$ where $\lambda'(v-1)=k'(k'-1)$.
\end{proof}

If equality holds in~\eqref{leqnec}, then $(v,k,\lambda)=(v,vk',v\lambda')$, where
$\lambda'(v-1)=k'(k'-1)$; thus, the parameters satisfy the numerical relation for a cube of propriety~$2$. In the special case $v=k$, we will see in Proposition~\ref{propls} that~$A$ is actually of propriety~$2$, but this is generally
not true. For example, only one of the $10$ cubes in the $(7,21,7)$-row of Table~\ref{tab4}
is of propriety~$2$. Two characteristic examples of these cubes are depicted in
Figure~\ref{fig5}.

Proposition~\ref{prop:1} gives a lower bound on~$\lambda$ for cubes of propriety~$3$:
\begin{equation}\label{lboundl}
\lambda \ge \left\lceil \frac{k(k-v)}{v(v-1)}\right\rceil.
\end{equation}
The inequality is strict if $k(k-v)$ is divisible by $v(v-1)$ and~$k$ is not divisible by~$v$, since equality would imply equality in Proposition~\ref{prop:1} and hence $v\mid k$.
The following is an improvement of Proposition~\ref{prop:1} if $k$ is not divisible by~$v$. Indeed, if $k=qv+r$ with $1\le r<v$, then the resulting lower bound on~$\lambda$ is larger by $r(v-r)/(v(v-1))$.

\begin{proposition}\label{prop:2}
If a cube $A\in \C^3_3(v,k,\lambda)$ exists and $k$ is not divisible by~$v$, then
%\begin{equation}\label{leqnec2}
%\lambda v(v-1)\geq v(k(\lfloor k/v\rfloor+\lceil k/v \rceil)-v\lfloor %k/v\rfloor\lceil k/v \rceil-k).
%\end{equation}
\begin{equation}\label{leqnec2}
\lambda (v-1)\geq k(\lfloor k/v\rfloor+\lceil k/v \rceil-1)-v\lfloor k/v\rfloor\lceil k/v \rceil.
\end{equation}
Equality holds if and only if every row-sum of~$A$ is either $\lfloor k/v\rfloor$ or $\lceil k/v\rceil$. In that case, the former sum occurs $\lceil k/v \rceil v-k$ times, and the latter sum
$k-\lfloor k/v\rfloor v$ times among parallel rows in any layer.
\end{proposition}

\begin{proof}
We use the same notation as in the proof of Proposition~\ref{prop:1}.
Since $k$ is not divisible by $v$, $\lfloor k/v\rfloor$ and $\lceil k/v \rceil$
are consecutive integers. Then, either $s_{xy}\leq \lfloor k/v\rfloor$ or
$s_{xy}\ge \lceil k/v \rceil$ holds, and thus
$(s_{xy}-\lfloor k/v\rfloor)(s_{xy}-\lceil k/v \rceil)\geq 0$ for all $x,y\in \{1,\ldots,v\}$.
This can be written as
\[
s_{xy}^2\ge \big(\lfloor k/v\rfloor + \lceil k/v \rceil\big)s_{xy} - \lfloor k/v\rfloor\,\lceil k/v \rceil.
\]
Under the constraint~\eqref{sumentries}, we have
\begin{align*}
\lambda \binom{v}{2} &= \sum_{x=1}^v \sum_{y=1}^v \binom{s_{xy}}{2}=\frac{1}{2}\left(\sum_{x=1}^v \sum_{y=1}^v s_{xy}^2-vk\right)\\
&\ge \frac{1}{2}\left(\sum_{x=1}^v \sum_{y=1}^v \Big( \big(\lfloor k/v\rfloor + \lceil k/v \rceil\big)s_{xy} - \lfloor k/v\rfloor\,\lceil k/v \rceil \Big) -vk\right)\\
&= \frac{1}{2}\Big(vk\big(\lfloor k/v\rfloor+\lceil k/v \rceil-1\big)-v^2\lfloor k/v\rfloor\lceil k/v \rceil\Big).
\end{align*}
This implies inequality~\eqref{leqnec2}, and equality holds if and only if
$(s_{xy}-\lfloor k/v\rfloor)(s_{xy}-\lceil k/v \rceil)=0$ for all $x,y\in \{1,\ldots,v\}$. The argument is also valid for row-sums in the direction of the $x$-axis and $y$-axis.

Assuming equality, let $n_1$ and $n_2$ be the number of occurrences of sums
$\lfloor k/v\rfloor$ and $\lceil k/v\rceil$ among the parallel rows in any layer,
respectively. Then, $n_1+n_2=v$ and $\lfloor k/v\rfloor n_1+\lceil k/v \rceil n_2=k$,
implying $n_1=\lceil k/v \rceil v-k$ and $n_2=k-\lfloor k/v\rfloor v$.
\end{proof}

For $(v,k)=(6,9)$, Proposition~\ref{prop:1} gives $\lambda\ge 1$, and Proposition~\ref{prop:2}
gives $\lambda\ge 2$. The following gives a nontrivial upper bound on $\lambda$ for $k<v$.

\begin{proposition}\label{prop:3}
If a cube $A\in \C^3_3(v,k,\lambda)$ exists and $k>1$, then
%\begin{equation*}%\label{leqnec3}
$\lambda (v-1)\leq k(k-1)-1$.
%\end{equation*}
\end{proposition}

\begin{proof}
Using the same notation as in the proof of Proposition~\ref{prop:1} and
$\sum_{y=1}^vs_{xy}=k$, we have
\[
\lambda \binom{v}{2} =\sum_{x=1}^v \sum_{y=1}^v \binom{s_{xy}}{2}=\frac{1}{2}\left(\sum_{x=1}^v \sum_{y=1}^v s_{xy}^2-vk\right)
\leq\frac{1}{2}\left(\sum_{x=1}^v \left(\sum_{y=1}^v s_{xy}\right)^2-vk\right)
= \frac{1}{2}\Big(v k^2-vk \Big).
\]
This implies $\lambda (v-1)\leq k(k-1)$, and equality holds only if $(s_{xy}/k)_{x,y=1}^v$
is a permutation matrix. After relabeling, we may assume $s_{xy}=k[x=y]$; then every row in an orthogonal direction contains at most one $1$-entry. Since the equality argument in that direction would require every nonzero row-sum to equal~$k$, simultaneous equality is possible only for $k=1$.
Therefore, we conclude that $\lambda (v-1)\leq k(k-1)-1$ holds.
\end{proof}

Let us now look at the extremal cases of the trivial inequality $0\le \lambda\le k$.

\begin{proposition}\label{proplameqk}
If $\C^3_3(v,k,\lambda)\neq\emptyset$ and $\lambda=k$, then $k=0$ or $k=v^2$.
\end{proposition}

\begin{proof}
Let us assume that $A\in \C^3_3(v,k,\lambda)$ is not the zero array,
that is, $A(x_0,y_0,z_0)=1$ for some $(x_0,y_0,z_0)$.
Because the total number of $1$-entries in a layer is equal to the
scalar product of parallel layers ($k=\lambda$), the $1$-entries
in parallel layers must appear in the same positions. Therefore,
the three rows through $(x_0,y_0,z_0)$ contain only $1$-entries,
i.e.\ $A(x,y_0,z_0)=A(x_0,y,z_0)=A(x_0,y_0,z)=1$ for all
$x,y,z\in\{1,\ldots,v\}$. From this it follows that all
entries of~$A$ are~$1$. The zero array has $k=0$, while
the all-$1$ array has $k=v^2$.
\end{proof}

\begin{proposition}\label{proplam0}
If $\C^3_3(v,k,\lambda)\neq\emptyset$ and $\lambda=0$, then $k\le v$.
\end{proposition}

\begin{proof}
This is a direct consequence of~\eqref{leqnec}.
\end{proof}

For $\lambda=0$, the next two propositions characterize the
extremal cases of $1\le k \le v$. Independently permuting
symbols $\{1,\ldots,v\}$ in the three coordinates of a
cube~$A$ is called \emph{isotopy}. Permuting the order of
the coordinates is called \emph{conjugation}.
Two cubes are considered \emph{equivalent} if one can be
transformed into the other by these operations;
see~\cite{KP25, KPT25, KR25} for more elaborate definitions.

\begin{proposition}\label{propidcube}
All cubes in $\C^3_3(v,1,0)$ are equivalent to $I(x,y,z)=[x=y=z]$.% $$I(x,y,z)=[x=y=z].$$
\end{proposition}

\begin{proof}
The $z=1$ layer of each $A\in \C^3_3(v,1,0)$ contains exactly one $1$-entry. By permuting symbols
in the $x$- and $y$-coordinates, we can put it in position $(1,1,1)$.
Next, the $z=2$ layer of~$A$ also contains exactly one $1$-entry.
This entry is neither in the $x=1$ layer, nor in the $y=1$ layer,
because they already contain the first $1$-entry. Therefore, we can
permute the symbols $\{2,\ldots,v\}$ in the $x$- and $y$-coordinates
to put the second $1$-entry in position $(2,2,2)$ without affecting
the first $1$-entry. Going on like this, we can put the remaining
$1$-entries in positions $(3,3,3)$ to $(v,v,v)$ and transform~$A$
to the ``identity cube''~$I$.
\end{proof}

\begin{figure}[!ht]
\ifthenelse{\boolean{images}}
{\ifthenelse{\boolean{lowres}}{\includegraphics[width=6cm]{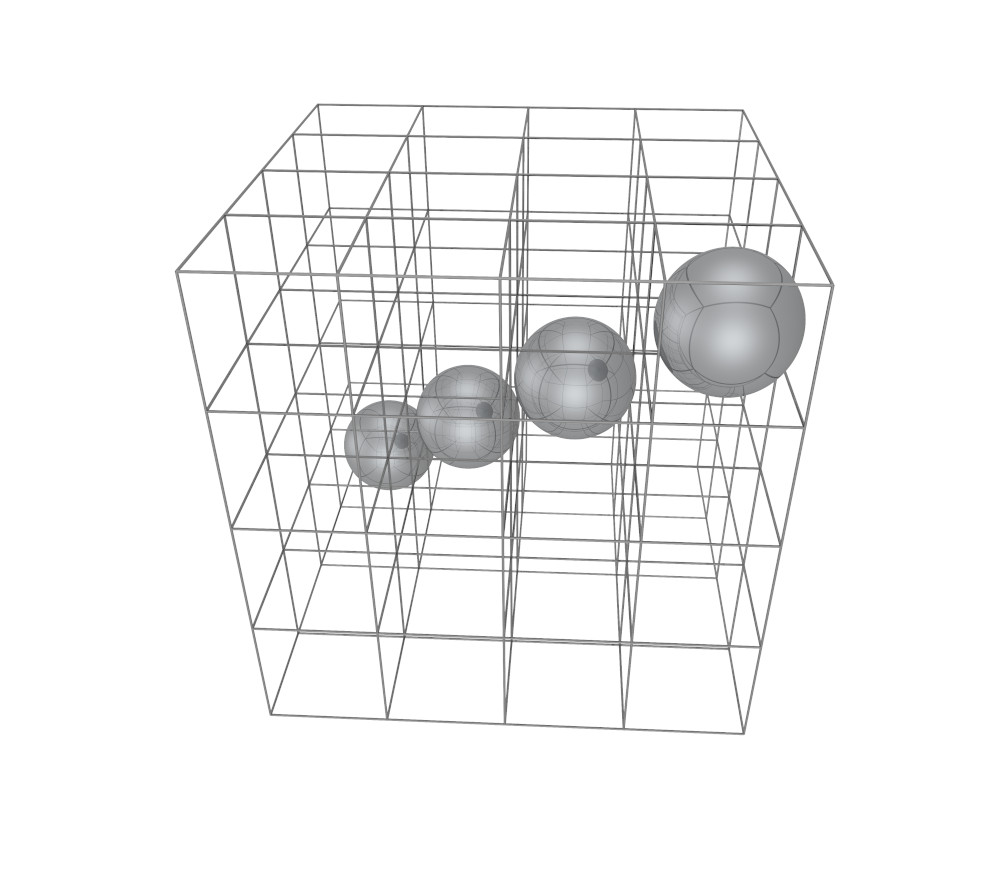}\\[-3mm]}{\includegraphics[width=6cm]{4-1-0.jpg}\\[-3mm]}
}{}
\caption{The identity cube $I\in \C^3_3(4,1,0)$.}\label{fig1}
\end{figure}

An image of the identity cube of order $v=4$ is shown in Figure~\ref{fig1}.
The spheres represent \hbox{$1$-entries}, and the empty cells represent $0$-entries.
The images in this paper were rendered using the ray tracing software
POV-Ray~\cite{POVRay}.

\begin{proposition}\label{propls}
All cubes in $\C^3_3(v,v,0)$ are of propriety $d=2$,
and they are equivalent to Latin squares of order~$v$.
\end{proposition}

\begin{proof}
Every layer of $A\in \C^3_3(v,v,0)$ contains $v$ entries~$1$.
Because of the $\lambda=0$ condition, they must appear in
different rows and columns. Hence, the layers of~$A$ are
permutation matrices, and $\C^3_3(v,v,0)=\C^3_2(v,1,0)$.
Cubes in $\C^3_2(v,1,0)$ are equivalent to Latin squares
of order~$v$ by the correspondence $A\leftrightarrow L$,
$A(x,y,z)=\left[ L(x,y)=z \right]$.
\end{proof}

\begin{figure}[t]
\ifthenelse{\boolean{images}}
{\ifthenelse{\boolean{lowres}}
{\includegraphics[width=5.6cm]{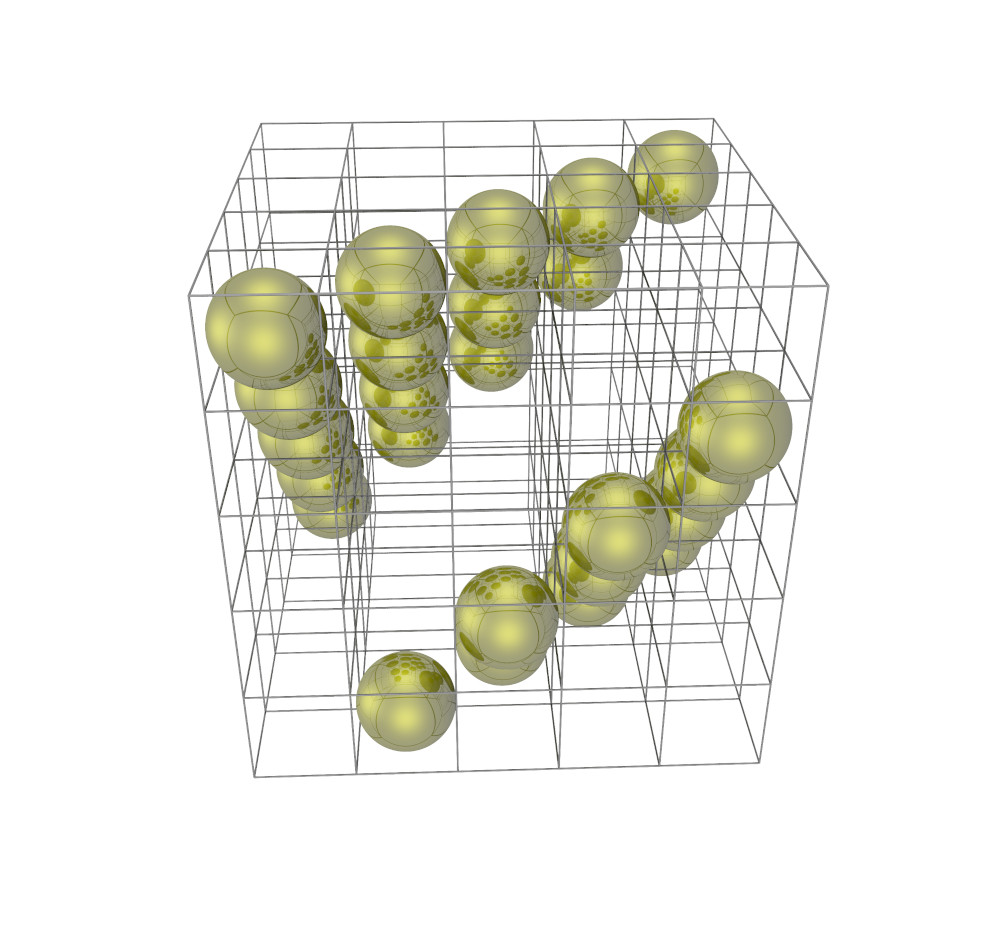} \includegraphics[width=5.6cm]{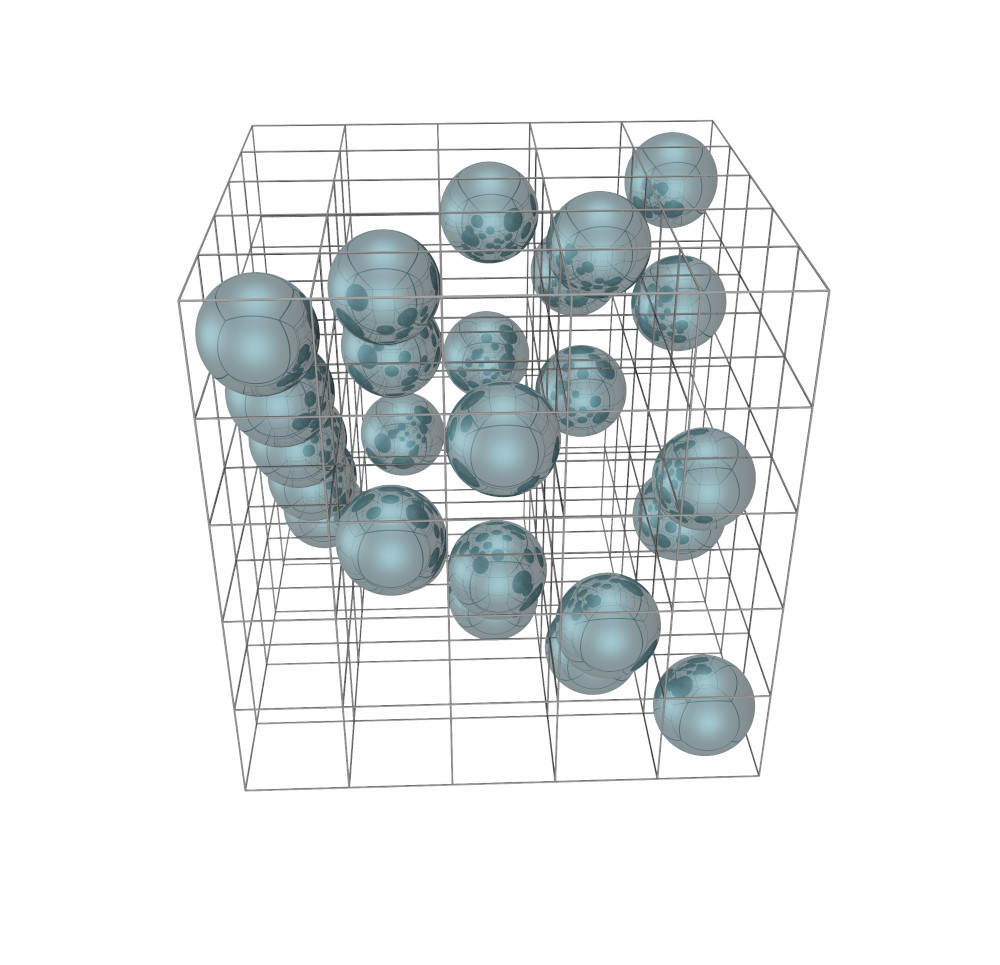}\\[-6mm]}
{\includegraphics[width=5.6cm]{5-5-0a.jpg} \includegraphics[width=5.6cm]{5-5-0b.jpg}\\[-6mm]}
}{}
$$\left[
\begin{array}{lllll}
 1 & 2 & 3 & 4 & 5 \\
 2 & 3 & 4 & 5 & 1 \\
 3 & 4 & 5 & 1 & 2 \\
 4 & 5 & 1 & 2 & 3 \\
 5 & 1 & 2 & 3 & 4
\end{array}
\right] \kern 20mm
\left[
\begin{array}{lllll}
 1 & 2 & 3 & 4 & 5 \\
 2 & 1 & 5 & 3 & 4 \\
 3 & 4 & 1 & 5 & 2 \\
 4 & 5 & 2 & 1 & 3 \\
 5 & 3 & 4 & 2 & 1
\end{array}
\right]$$
\caption{Two nonequivalent $\C^3_3(5,5,0)$-cubes and the corresponding Latin squares.}\label{fig2}
\end{figure}

\begin{table}[!b]
\vskip 7mm
\begin{tabular}{cccc}
\hline
$v$ & $k$ & $\lambda$ & No.\ cubes\\
\hline
3 & 2 & 0 & 2\\
3 & 3 & 0 & 1\\
 & & 1 & 1\\
3 & 4 & 1 & 1\\
\hline
4 & 2 & 0 & 8\\
4 & 3 & 0 & 4\\
4 & 4 & 0 & 2\\
4 & 5 & $\lambda$ & 0\\
4 & 6 & 2 & 117\\
4 & 7 & 2 & 4\\
 & & 3 & 2\\
4 & 8 & 4 & 19\\
\hline
5 & 2 & 0 & 23\\
5 & 3 & 0 & 251\\
5 & 4 & 0 & 40\\
\hline
6 & 2 & 0 & 157\\
\hline
\end{tabular}
\vskip 5mm
\caption{Numbers of nonequivalent $\C^3_3(v,k,\lambda)$-cubes.}\label{tab1}
\end{table}

The correspondence is illustrated in Figure~\ref{fig2} for $v=5$.
Main classes of Latin squares correspond to equivalence classes
of $\C^3_3(v,v,0)$-cubes; see~\cite{KD15} or~\cite{LM98} for
definitions related to Latin squares. Main classes of Latin
squares were enumerated by computer calculations up to
$v=11$~\cite{HKO11, MMM07}.

We performed a computer enumeration of $\C^3_3(v,k,\lambda)$-cubes
for $v\le 4$; $v=5$, $k\le 4$; and $v=6$, $k=2$. For fixed~$v$
and~$k$, we essentially generate all $v\times v\times v$ binary
arrays with~$k$ entries~$1$ in every layer. We then discard the
arrays not satisfying the $\lambda$-condition, and eliminate
equivalent cubes by applying the \texttt{CubeFilter} command from
the GAP~\cite{GAP} package PAG~\cite{PAG}. The command transforms
incidence cubes to equivalent colored graphs and calls
nauty/Traces~\cite{MP14} to produce canonical labelings;
see~\cite[Section~2]{KPT25} for implementation details.
To speed up the calculation, isomorph rejection is also performed
on partial cubes, consisting of one or two layers.

Results of the enumeration are presented in Table~\ref{tab1}.
We only consider parameters satisfying $0\le k \le v^2/2$,
because the remaining cases are covered by taking complements.
For each such pair $(v,k)$, all possible values of~$\lambda$
are given in consecutive rows of Table~\ref{tab1}. Two
values of~$\lambda$ occur for $(3,3)$ and $(4,7)$, while
for $(4,5)$ there are no possible values of~$\lambda$.
An online version of the table, with links to files
containing the $\C^3_3(v,k,\lambda)$-cubes in GAP-readable
format, is available on the web page %$\ldots$
\begin{center}
\url{https://web.math.pmf.unizg.hr/~krcko/results/cubes3.html}.
\end{center}

The bound~\eqref{lboundl} is tight for all values $(v,k)$ from Table~\ref{tab1} except
$(4,8)$. For $(v,k)=(4,8)$, the bound gives $\lambda\ge 3$, but cubes exist only
for $\lambda=4$. The set $\C^3_3(4,8,4)$ is closed under complementation. Nine cubes
are equivalent to their complements, while $10$ cubes form nonequivalent complementary pairs,
accounting for a total of $19$ nonequivalent $\C^3_3(4,8,4)$-cubes.

\section{Difference sets}\label{sec3}

Let $G$ be a group of order~$v$. We will use additive notation,
but~$G$ is not assumed to be commutative.
In this section, cubes are regarded as functions $A:G^3\to \{0,1\}$.
The \emph{full autotopy group} of~$A$ contains all isotopies
from~$A$ to~$A$: $\Atop(A)=\{(\alpha,\beta,\gamma)\in S(G)^3 \mid
A(\alpha(x),\beta(y),\gamma(z))=A(x,y,z),\, \forall x,y,z\in G\}$.
Here, $S(G)$ denotes the set of all permutations of the elements
of~$G$. The group operation in $\Atop(A)$ is coordinatewise
composition. The phrase ``an autotopy group of~$A$'' refers
to a subgroup of~$\Atop(A)$.

A \emph{$(v,k,\lambda)$ difference set} in~$G$ is a subset
$D\subseteq G$ of size $|D|=k$ such that, for every
$g\in G\setminus\{0\}$, there are exactly~$\lambda$
choices of $x,x'\in D$ with $x-x'=g$. Then, the \emph{development}
$\dev D=\{D+g \mid g\in G\}$ comprises blocks of
a symmetric $(v,k,\lambda)$ design. Symmetric designs constructed
from difference sets in~$G$ are characterized by having an automorphism
group isomorphic to~$G$ acting sharply transitively on the points
and blocks. A \emph{multiplier} of~$D$ is a group automorphism
$\varphi \in \Aut(G)$ such that $\varphi(D)=D+a$ for some $a\in G$.
The set of all multipliers $\Mult(D)$ is a subgroup of $\Aut(G)$,
and each multiplier induces an automorphism of $\dev D$.
See \cite[Chapter~VI]{BJL99}, \cite[Chapter~9]{IS06}, or
\cite[Chapter~4]{ESL83} for more results about difference sets.

By \cite[Theorem~3.1]{KPT25}, the following formula involving the
Iverson symbol defines a $\C^3_2(v,k,\lambda)$-cube of propriety~$2$
from a $(v,k,\lambda)$ difference set~$D$:
\begin{equation}\label{ds2}
A(x,y,z)=[x+y+z\in D].
\end{equation}
Examples constructed in this way are called \emph{difference
cubes}. By~\cite[Theorem~3.4]{KPT25}, difference cubes of
dimension~$3$ have an autotopy group isomorphic to a semidirect
product $G^2\rtimes \Mult(D)$.

A variant of difference sets for $n$-dimensional projection cubes
was introduced in \cite[Definition~3.1]{KR25}. Using a suitable
normalization, the three-dimensional case can be defined as a subset
of pairs $D\subseteq G^2$ of size $|D|=k$ such that the following
subsets of~$G$ are ``ordinary'' $(v,k,\lambda)$ difference
sets in~$G$:
\begin{enumerate}[label=\arabic*.]
\item the first coordinates $\{x\in G \mid \exists (x,y)\in D\}$,
\item the second coordinates $\{y\in G \mid \exists (x,y)\in D\}$, and
\item differences of the coordinates $\{x-y \mid (x,y)\in D\}$.
\end{enumerate}
Then, by \cite[Proposition~3.2]{KR25}, the development
\begin{equation}\label{dsdev}
\dev D := \{ (x+g,y+g,g) \mid (x,y)\in D,\, g\in G\}
\end{equation}
is the support of a $\P^3(v,k,\lambda)$-cube.
Here and in the sequel, the \emph{support} of a function
$A:G^3\to \{0,1\}$ means the set of triples
$\supp A:=\{(x,y,z)\in G^3 \mid A(x,y,z)=1\}$.

We now define the corresponding concept for three-dimensional
symmetric designs of propriety~$3$.

\begin{definition}\label{ds3def}
A subset of pairs $D\subseteq G^2$ of size $|D|=k$ is a \emph{$(v,k,\lambda)$
difference set of propriety~$3$} if, for every $g\in G\setminus\{0\}$,
the following hold.
\begin{enumerate}[label=\textup{\arabic*.}]
\item There are exactly $\lambda$ choices of $(x,y), (x',y')\in D$
with $-x+x' = -y+y' = g$.
\item There are exactly $\lambda$ choices of $(x,y), (x',y')\in D$
with \hbox{$x-x'=g$} and $y=y'$.
\item There are exactly $\lambda$ choices of $(x,y), (x',y')\in D$
with $x=x'$ and $y-y'=g$.
\end{enumerate}
\end{definition}

Note that condition~1 involves the left differences $-x+x'$,
$-y+y'$, while conditions~2 and~3 involve the right
differences $x-x'$ and $y-y'$. This would have to be reversed
if left translation was used in~\eqref{dsdev} instead
of right translation; cf.~\cite[Example~3.4]{KR25}.

\begin{theorem}\label{ds3tm}
If $D\subseteq G^2$ satisfies the conditions of Definition~\ref{ds3def},
then the development~\eqref{dsdev} is the support of a
$\C^3_3(v,k,\lambda)$-cube.
\end{theorem}

\begin{proof}
Let $A:G^3\to\{0,1\}$ be the array with $\supp A=\dev D$.
The number of $1$-entries in a layer is the number of triples
in~\eqref{dsdev} with a fixed coordinate. If the third
coordinate is fixed, this is obviously $|D|=k$. If the first
coordinate is fixed, e.g.\ $x+g=c$, then $g=-x+c$ and we are
counting pairs $(y-x+c,-x+c)$ with $(x,y)\in D$. The number
of such pairs is also~$k$, because $(x,y)\mapsto (y-x+c,-x+c)$
is an injective mapping. The argument is similar if the second
coordinate is fixed.

To compute the scalar product of two parallel layers, we fix
one coordinate to two values $c_1\neq c_2$ and count how many
times pairs of the remaining coordinates in~\eqref{dsdev}
coincide. If the third coordinate is fixed, we are counting
pairs $(x,y), (x',y')\in D$ with $(x+c_1,y+c_1)=(x'+c_2,y'+c_2)$.
This is equivalent to $-x+x'=-y+y'=c_1-c_2$, and the number of
pairs is~$\lambda$ by condition~1. If the first coordinate
is fixed, we are counting pairs $(x,y), (x',y')\in D$ with
$(y-x+c_1,-x+c_1)=(y'-x'+c_2,-x'+c_2)$. This is equivalent to
$x-x'=c_1-c_2$, $y=y'$, and the result is~$\lambda$ by
condition~2. Condition~3 is invoked if the second coordinate
is fixed. Thus, the scalar product of any pair of parallel
layers is~$\lambda$, and $A\in \C^3_3(v,k,\lambda)$.
\end{proof}

Difference sets for three-dimensional projection cubes are a
special case of Definition~\ref{ds3def}.

\begin{proposition}
Let $D\subseteq G^2$ be a difference set for
$\P^3(v,k,\lambda)$-cubes. Then, $D$ is also a
$(v,k,0)$ difference set of propriety~$3$.
\end{proposition}

\begin{proof}
By the definition of difference sets for projection $3$-cubes,
the mappings $(x,y)\mapsto x$, $(x,y)\mapsto y$,
and $(x,y)\mapsto x-y$ take the pairs in~$D$ to three
$(v,k,\lambda)$ difference sets in~$G$. The domain and the
images are sets of size~$k$, therefore the three mappings
are injective. Assume there are two pairs $(x,y),(x',y')\in D$
with $-x+x'=-y+y'=g\neq 0$. Then, $x-y=x'-y'$ contradicts the
injectivity of $(x,y)\mapsto x-y$. Therefore, condition~1
of Definition~\ref{ds3def} holds for $\lambda=0$. Similarly,
injectivity of $(x,y)\mapsto y$ implies condition~2 for
$\lambda=0$, and injectivity of $(x,y)\mapsto x$ implies
condition~3 for $\lambda=0$.
\end{proof}

\begin{figure}[!b]
\ifthenelse{\boolean{images}}
{\ifthenelse{\boolean{lowres}}
{\includegraphics[width=6cm]{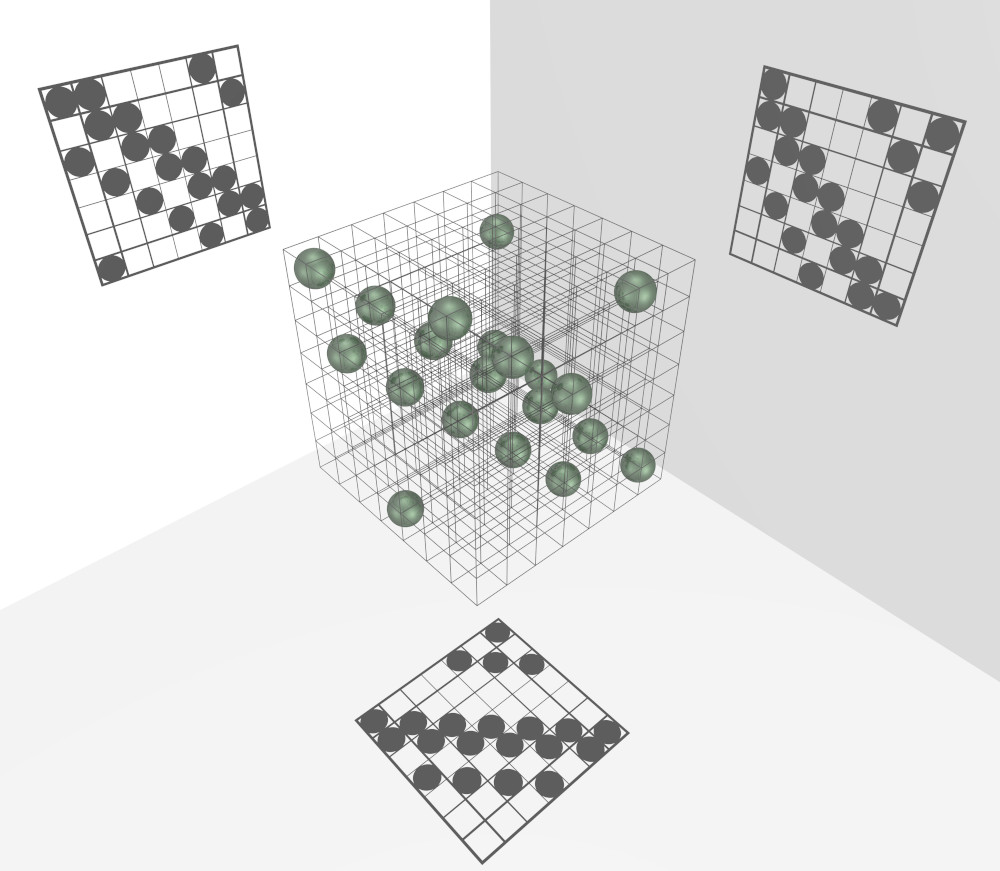}\hskip 6mm
\includegraphics[width=6cm]{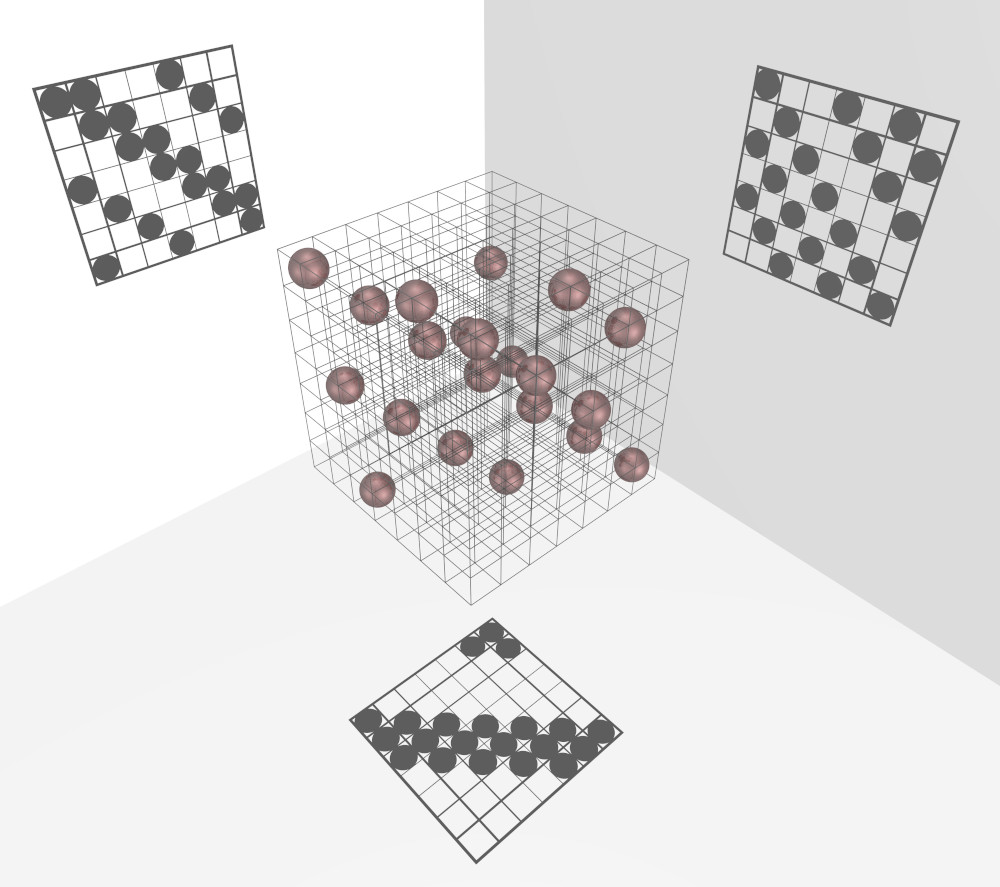}}
{\includegraphics[width=6cm]{7-3-0d1.jpg}\hskip 6mm
\includegraphics[width=6cm]{7-3-0d2.jpg}}
}{}
\caption{Two cubes in $\C^3_3(7,3,0)$.}\label{fig4}
\end{figure}

\begin{example}
Let $G=\Z_7$ be the integers with addition modulo~$7$ and
\[
D_1 = \big\{ (0, 0),\, (1, 2),\, (3, 6) \big\},\quad
D_2 = \big\{ (0, 0),\, (1, 3),\, (2, 6) \big\}.
\]
The sets $D_1$ and $D_2$ are both $(7,3,0)$ difference
sets of propriety~$3$. The set $D_1$ is also a difference
set for $\P^3(7,3,1)$-cubes, while $D_2$ is not.
\end{example}

The cubes obtained from $D_1$ and $D_2$ are shown in
Figure~\ref{fig4}. Three light sources are placed along
the coordinate axes, making the projections visible as
shadows. In the left image, the shadows are incidence
matrices of the Fano plane, while in the right image
they are not.

\begin{example}
Let $G=\Z_7$ be the integers with addition modulo~$7$ and
\begin{align*}
D_3 = \big\{\, & ( 0, 1 ),\, ( 0, 2 ),\, ( 0, 4 ),\, ( 1, 0 ),\, ( 1, 1 ),\, ( 1, 3 ),\, ( 2, 0 ),\, ( 2, 2 ),\, ( 2, 6 ),\, ( 3, 1 ),\, ( 3, 5 ),\,\\[0mm]
         &  ( 3, 6 ),\, ( 4, 0 ),\, ( 4, 4 ),\, ( 4, 5 ),\, ( 5, 3 ),\, ( 5, 4 ),\, ( 5, 6 ),\, ( 6, 2 ),\, ( 6, 3 ),\, ( 6, 5 )  \,\big\},\\[2mm]
D_4 = \big\{\, & ( 0, 1 ),\, ( 0, 2 ),\, ( 0, 4 ),\, ( 1, 0 ),\, ( 1, 2 ),\, ( 1, 4 ),\, ( 2, 0 ),\, ( 2, 1 ),\, ( 2, 4 ),\, ( 3, 3 ),\, ( 3, 5 ),\,\\[0mm]
             &  ( 3, 6 ),\, ( 4, 0 ),\, ( 4, 1 ),\, ( 4, 2 ),\, ( 5, 3 ),\, ( 5, 5 ),\, ( 5, 6 ),\, ( 6, 3 ),\, ( 6, 5 ),\, ( 6, 6 )\, \big\}.
\end{align*}
The sets $D_3$ and $D_4$ are both $(7,21,7)$ difference sets of propriety~$3$.
\end{example}

\begin{figure}[!ht]
\ifthenelse{\boolean{images}}
{\ifthenelse{\boolean{lowres}}
{\includegraphics[width=6cm]{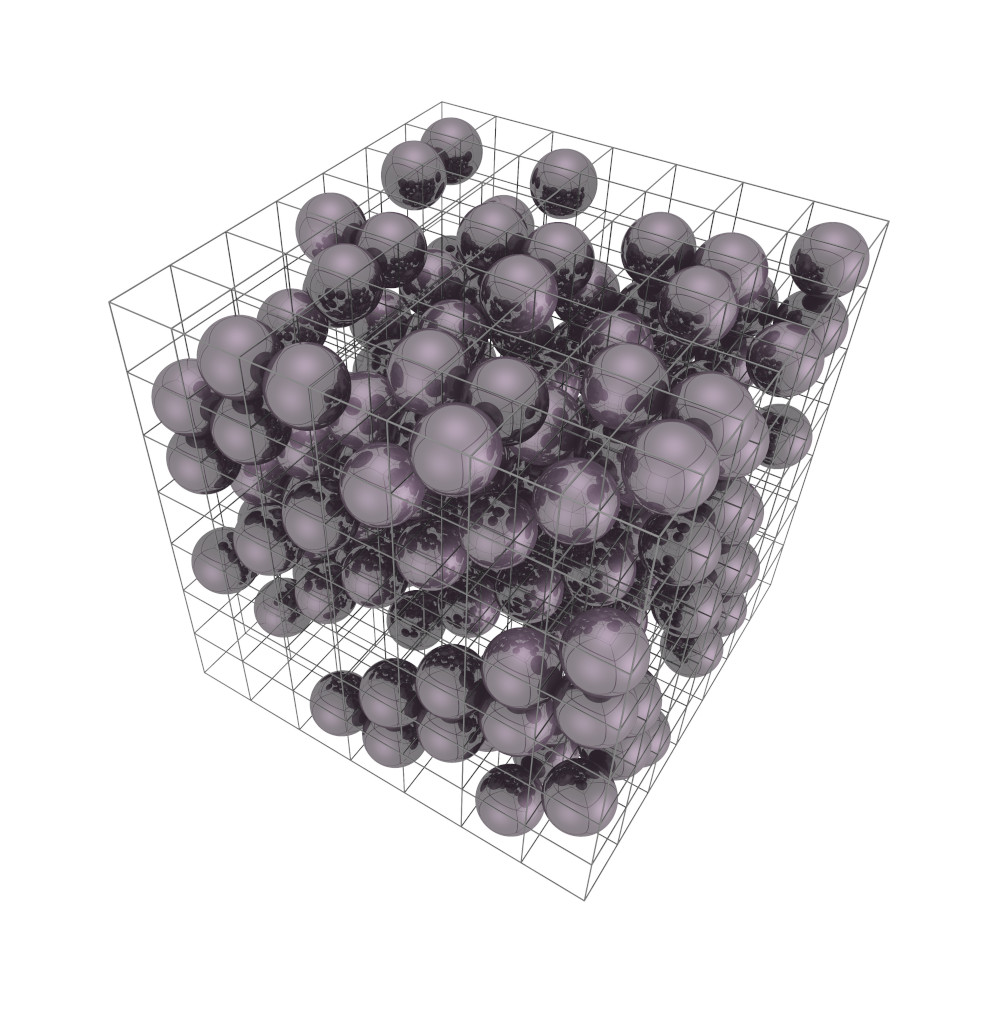}\includegraphics[width=6cm]{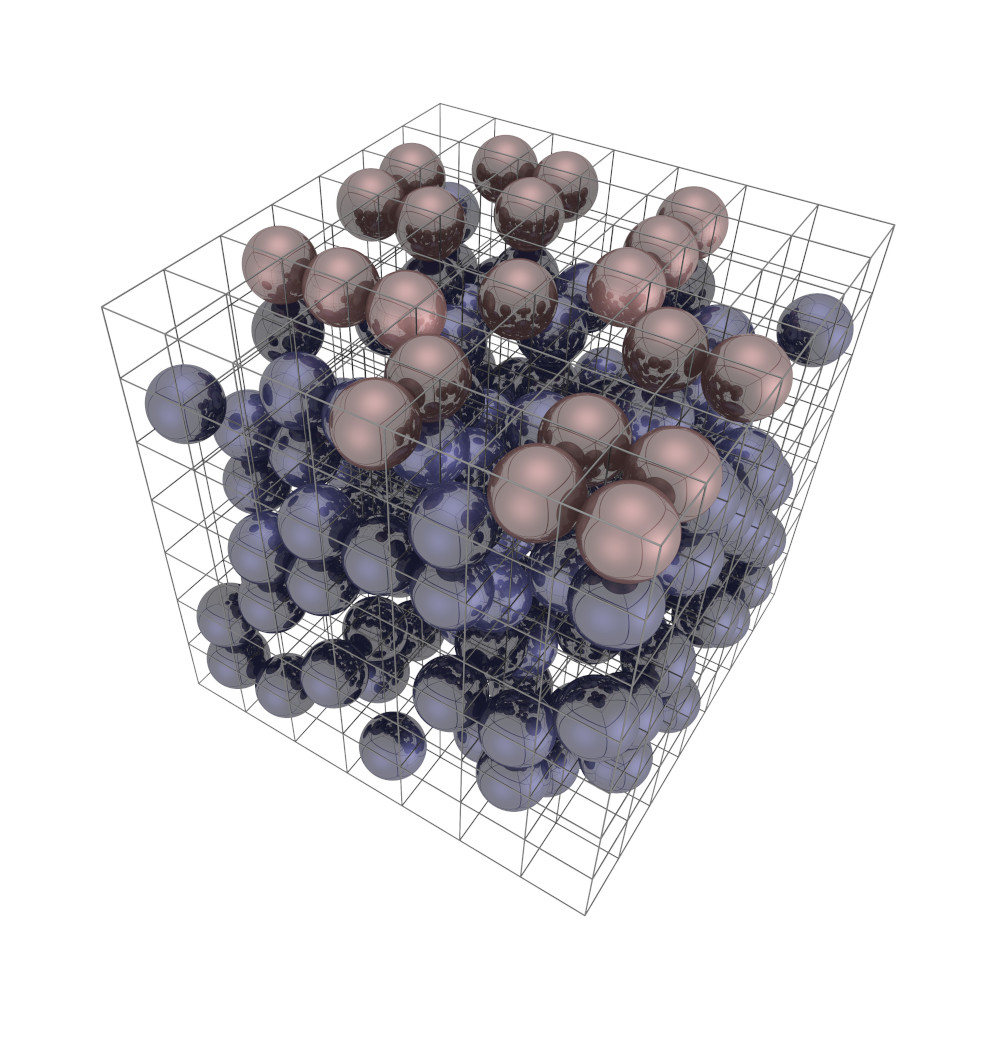}}
{\includegraphics[width=6cm]{7-21-7d3.jpg}\includegraphics[width=6cm]{7-21-7d4.jpg}}
}{}
\vskip -3mm
\caption{Two cubes in $\C^3_3(7,21,7)$.}\label{fig5}
\end{figure}

The development $\dev D_3$ actually supports a $\C^3_2(7,3,1)$-cube
of propriety~$2$. The development $\dev D_4$ supports a
$\C^3_3(7,21,7)$-cube which is not of propriety~$2$. The cubes are
shown in Figure~\ref{fig5}; the top horizontal layer in the right
image is clearly not an incidence matrix of the Fano plane.

The cube obtained from~$D_3$ is equivalent to a difference
cube~\eqref{ds2} constructed from an ordinary $(7,3,1)$ difference set.
Surprisingly, the difference cube construction~\eqref{ds2} is generally
not a special case of Theorem~\ref{ds3tm}. To give counterexamples,
we need to characterize $\C^3_3(v,k,\lambda)$-cubes coming from
difference sets of propriety~$3$.

\begin{theorem}\label{ds3char}
A cube $A\in \C^3_3(v,k,\lambda)$ can be constructed from a $(v,k,\lambda)$
difference set of propriety~$3$ in~$G$ if and only if $\Atop(A)$ has
a subgroup isomorphic to~$G$ acting sharply transitively on the three
coordinates.
\end{theorem}

\begin{proof}
The proof is essentially the same as proofs of Propositions~3.5
and~3.7 in~\cite{KR25}. If $A$ is constructed from a difference set~$D$ of propriety~$3$, then the right-regular action of~$G$, acting diagonally on the three coordinates, gives a subgroup of $\Atop(A)$ isomorphic to~$G$ and acting sharply transitively on all three coordinates.
For the converse, suppose that~$A$ has an autotopy group~$G$ acting sharply
transitively on the three coordinates.
If the indexing set is $\{1,\ldots,v\}$, we can enumerate the elements
of $G$ as $g_i:=(\alpha_i,\beta_i,\gamma_i)$, where the components are
permutations in $S_v$ with $\alpha_i(1)=i$, $\beta_i(1)=i$, and
$\gamma_i(1)=i$, for $i=1,\ldots,v$. This follows from sharp
transitivity and by applying isotopy to~$A$. We can now identify
$i\leftrightarrow g_i$ so that $A:G^3\to \{0,1\}$ and $G$ acts by
right translation. The support of~$A$ contains $vk$ triples
and splits into $k$ orbits of size~$v$. From each orbit we take a
representative of the form $(x_i,y_i,0)$, $i=1,\ldots,k$,
where $x_i,y_i\in G$ and $0\in G$ is the neutral element.
Then, $D:=\{(x_i,y_i) \mid i=1,\ldots,k\}$ satisfies
$\dev D=\supp A$, and this implies the conditions
of Definition~\ref{ds3def}.
\end{proof}

We use multiplicative notation in the next example.

\begin{example}\label{exQ16}
Let $G=\langle a,b\mid a^8=1,\, b^2=a^4,\, ba=a^7b\rangle$
be the generalized quaternion group of order~$16$.
Then, $D:=\{1,a^6,b,b^3,ab^3,a^3b^3\}$ is a $(16,6,2)$ difference
set in~$G$. Let $A(x,y,z)=[x\cdot y\cdot z\in D]$ be the
corresponding difference cube in $\C^3_2(16,6,2)$. Using
the package PAG~\cite{PAG}, we computed $\Atop(A)$ and found
that it is a group of order $1024$ isomorphic to
$(G\times G)\rtimes (C_2\times C_2)$. More
computations in GAP showed that $\Atop(A)$ has exactly $1557$
subgroups of order~$16$ up to conjugation, none of them
acting transitively on all three coordinates. Therefore,
by Theorem~\ref{ds3char}, the difference cube $A$ does not
arise from a $(16,96,32)$ difference set of propriety~$3$.
\end{example}

By~\cite[Table~1 and Proposition~5.1]{KPT25}, there are
exactly~$27$ nonequivalent difference cubes in $\C^3_2(16,6,2)$.
We performed computations in GAP as in Example~\ref{exQ16}
and established that~$8$ of these cubes can be constructed from
difference sets of propriety~$3$, while $19$ do
not arise in this way.

\begin{table}[!b]
\begin{tabular}{cccc}
\hline
$v$ & $k$ & $\lambda$ & \rule{4mm}{0mm}$C_v$\rule{4mm}{0mm}\\
\hline
3 & 2 & 0 & 1\\
3 & 3 & 0 & 1\\
 & & 1 & 1\\
3 & 4 & 1 & 1\\
\hline
5 & 2 & 0 & 1\\
5 & 3 & 0 & 2\\
5 & 4 & 0 & 1\\
5 & 5 & 0 & 1 \\
 & & 1 & 1\\
5 & 6 & 1 & 4\\
5 & 7 & 1 & 1\\
 & & 2 & 6\\
5 & 8 & 2 & 2\\
5 & 9 & 2 & 2\\
 & & 3 & 2\\
 & & 4 & 1\\
5 & 10 & 3 & 6\\
 & & 4 & 15\\
5 & 11 & 4 & 14\\
 & & 5 & 7\\
5 & 12 & 5 & 5 \\
 & & 6 & 1\\
\hline
\end{tabular}\hskip 10mm
\begin{tabular}{cccc}
\hline
$v$ & $k$ & $\lambda$ & \rule{4mm}{0mm}$C_v$\rule{4mm}{0mm}\\
\hline
7 & 2 & 0 & 2\\
7 & 3 & 0 & 7\\
7 & 4 & 0 & 10\\
7 & 5 & 0 & 9\\
7 & 6 & 0 & 4\\
 & & 1 & 5\\
7 & 7 & 0 & 2\\
 & & 1 & 26\\
7 & 8 & 1 & 104\\
 & & 2 & 1\\
7 & 9 & 1 & 158\\
 & & 2 & 11\\
7 & 10 & 1 & 24\\
 & & 2 & 152\\
 & & 3 & 1\\
7 & 11 & 2 & 1258\\
 & & 3 & 11\\
7 & 12 & 2 & 659\\
 & & 3 & 542\\
 & & 4 & 1\\
7 & 13 & 2 & 5\\
 & & 3 & 3875\\
 & & 4 & 56\\
7 & 14 & 3 & 1943\\
 & & 4 & 2021\\
 & & 5 & 17\\
\hline
\end{tabular}
\vskip 5mm
\caption{Numbers of nonequivalent $\C^3_3(v,k,\lambda)$-cubes coming
from cyclic difference sets.}\label{tab2}
\end{table}

\begin{table}[t]
\begin{tabular}{cccc@{}c@{\rule{2mm}{0mm}}c}
\hline
$v$ & $k$ & $\lambda$ & \rule{3mm}{0mm}$C_4$\rule{3mm}{0mm} & $C_2\times C_2$ & Total\\
\hline
4 & 2 & 0 & 1 & 1 & 2\\
4 & 3 & 0 & 1 & 1 & 2\\
4 & 4 & 0 & 0 & 1 & 1\\
4 & 5 & $\lambda$ & 0 & 0 & 0\\
4 & 6 & 2 & 5 & 6 & 7\\
4 & 7 & 2 & 1 & 1 & 1\\
 & & 3 & 0 & 0 & 0\\
4 & 8 & 4 & 0 & 1 & 1\\
\hline
\end{tabular}\hskip 5mm
\begin{tabular}{cccccc}
\hline
$v$ & $k$ & $\lambda$ & \rule{3mm}{0mm}$C_6$\rule{3mm}{0mm} & \rule{1mm}{0mm}$S_3$\rule{1mm}{0mm} & Total\\
\hline
6 & 2 & 0 & 3 & 2 & 4\\
6 & 3 & 0 & 7 & 4 & 9\\
6 & 4 & 0 & 5 & 4 & 9\\
6 & 5 & 0 & 2 & 1 & 3\\
6 & 6 & $\lambda$ & 0 & 0 & 0\\
6 & 7 & $\lambda$ & 0 & 0 & 0\\
6 & 8 & 2 & 10 & 3 & 13\\
6 & 9 & 2 & 158 & 42 & 200\\
6 & 10 & 2 & 71 & 12 & 83\\
6 & 11 & 2 & 2 & 0 & 2\\
 & & 4 & 3 & 0 & 3\\
6 & 12 & 4 & 283 & 42 & 322\\
6 & 13 & 4 & 91 & 13 & 99\\
6 & 14 & 4 & 9 & 1 & 9\\
 & & 6 & 58 & 10 & 66\\
6 & 15 & 6 & 147 & 29 & 160\\
6 & 16 & 6 & 109 & 7 & 116\\
 & & 8 & 10 & 6 & 14\\
6 & 17 & 8 & 180 & 17 & 193\\
6 & 18 & 8 & 264 & 30 & 292\\
 & & 10 & 4 & 2 & 4\\
\hline
\end{tabular}
\vskip 5mm
\caption{Numbers of nonequivalent $\C^3_3(v,k,\lambda)$-cubes coming
from difference sets.}\label{tab3}
\end{table}

Cubes coming from difference sets are more easily enumerated than
general $\C_3^3(v,k,\lambda)$-cubes. The results of a computer
enumeration up to equivalence are summarized in
Table~\ref{tab2} for $v=3,5,7$, and in Table~\ref{tab3} for $v=4,6$.
The former orders are prime and there is only the cyclic
group~$C_v$, while for the latter orders there
is one other group. For example, in the $(v,k,\lambda)=(4,6,2)$ row
of Table~\ref{tab3}, there are $5$ nonequivalent cubes coming
from difference sets of propriety~$3$ in the cyclic group~$C_4$,
$6$~in the Klein four-group $C_2\times C_2$, and~$7$ cubes in total.
This means that $5+6-7=4$ cubes in $\C_3^3(4,6,2)$ can be constructed
from difference sets of propriety~$3$ in both groups. Online versions
of the tables with links to difference sets and cubes in GAP-readable
format are available on our web page.

To perform the enumeration, given~$G$ and~$k$, we essentially
generate all $k$-subsets of~$G^2$ and check which are difference
sets of propriety~$3$. Thus, all possible $\lambda$-values are
obtained. For~$G=C_7$, a complete enumeration was possible only
for $k\le 14$. We now introduce multipliers to facilitate
enumeration for $k>14$.

\begin{lemma}\label{lmtrans}
Let $D\subseteq G^2$ be a $(v,k,\lambda)$ difference set of
propriety~$3$. Then, for each $(a,b)\in G^2$, the set
$(a,b)+D:=\{(a+x,b+y) \mid (x,y)\in D\}$
is also a $(v,k,\lambda)$ difference set of propriety~$3$.
The cubes supported by $\dev D$ and $\dev((a,b)+D)$ are
isotopic.
\end{lemma}

\begin{proof}
The triple $(\alpha,\beta,\gamma)\in S(G)^3$ with components
$\alpha(x) = a+x$, $\beta(y)=b+y$, $\gamma(z)=z$, applied coordinatewise,
maps $\dev D$ to $\dev((a,b)+D)$. From this it follows that $(a,b)+D$ is
also a $(v,k,\lambda)$ difference set of propriety~$3$, and $(\alpha,\beta,\gamma)$
is an isotopy of the supported cubes.
\end{proof}
The following is an analogue of multipliers for difference sets of propriety~$3$.

\begin{definition}
Let $D\subseteq G^2$ be a $(v,k,\lambda)$ difference set of
propriety~$3$. A \emph{multiplier} of~$D$ is a group
automorphism $\varphi \in \Aut(G)$ such that
$\varphi(D)=(a,b)+D$ for some $(a,b)\in G^2$.
\end{definition}

\begin{proposition}
A multiplier of $D\subseteq G^2$ induces an autotopy of the
cube supported by $\dev D$. The set of all multipliers
$\Mult(D)$ is a subgroup of $\Aut(G)$.
\end{proposition}

\begin{proof}
Let $\varphi$ be a multiplier of~$D$ with $\varphi(D)=(a,b)+D$.
Similarly as in Lemma~\ref{lmtrans}, the functions
$\alpha(x)=-a+\varphi(x)$, $\beta(y)=-b+\varphi(y)$,
and $\gamma(z)=\varphi(z)$ are permutations of~$G$.
Applied coordinatewise, $(\alpha,\beta,\gamma)$
maps $\dev D$ to itself and is therefore an
autotopy of the supported cube induced by~$\varphi$.

Given two multipliers $\varphi_i(D)=(a_i,b_i)+D$, $i=1,2$, it
is readily checked that $(\varphi_1\circ \varphi_2^{-1})(D)=
(a,b)+D$ for $a=-(\varphi_1\circ \varphi_2^{-1})(a_2)+a_1$
and $b=-(\varphi_1\circ \varphi_2^{-1})(b_2)+b_1$. Therefore,
$\varphi_1\circ \varphi_2^{-1}$ is also a multiplier
of~$D$ and $\Mult(D)$ is a subgroup of $\Aut(G)$.
\end{proof}

An important property of a multiplier $\varphi$ of an ordinary
difference set~$D$ is that there exists a
translate $D+g$ fixed by~$\varphi$
\cite[Lemma~VI.2.4]{BJL99}. Regarding difference sets of
propriety~$3$, we use an additional assumption on~$\varphi$
to prove the analogous property.

\begin{lemma}
Let $D\subseteq G^2$ be a $(v,k,\lambda)$ difference set of
propriety~$3$ and $\varphi \in \Aut(G)$ a multiplier of~$D$.
If the neutral element $0$ is the only element of~$G$ fixed
by~$\varphi$, then there exists a pair $(g,h)\in G^2$ such
that
\begin{equation}\label{eqfixed}
\varphi((g,h)+D)=(g,h)+D.
\end{equation}
\end{lemma}

\begin{proof}
Let $\varphi(D)=(a,b)+D$ and define $\psi:G\to G$ by
$\psi(x)=-x+\varphi(x)$. If $\psi(x)=\psi(y)$, then $-x+\varphi(x)=-y+\varphi(y)$ implies $y-x=\varphi(y-x)$, so $x=y$ by the assumption on~$\varphi$. Thus $\psi$ is injective, and hence bijective because~$G$ is finite. Therefore, there
exist $g,h\in G$ such that $\psi(g)=-a$ and $\psi(h)=-b$,
i.e.\ $\varphi(g)=g-a$ and $\varphi(h)=h-b$. It is now
easily verified that~\eqref{eqfixed} holds.
\end{proof}

\begin{table}[!b]
\begin{tabular}{cccccc}
\hline
$v$ & $k$ & $\lambda$ & $C_7\rtimes C_2$ & $C_7\rtimes C_3$ & Total\\
\hline
7 & 15 & 3 & 0 & 5 & 5 \\
  &    & 4 & 0 & 63 & 63 \\
  &    & 6 & 0 & 2 & 2 \\
7 & 16 & 4 & 0 & 58 & 58 \\
  &    & 5 & 0 & 4 & 4 \\
  &    & 6 & 1 & 1 & 2 \\
7 & 17 & 6 & 1 & 0 & 1 \\
7 & 18 & 5 & 0 & 24 & 24 \\
  &    & 6 & 3 & 5 & 5 \\
  &    & 7 & 1 & 36 & 36 \\
7 & 19 & 6 & 3 & 82 & 82 \\
  &    & 7 & 0 & 3 & 3 \\
  &    & 8 & 0 & 2 & 2 \\
  &    & 9 & 1 & 1 & 1 \\
7 & 20 & $\lambda$ & 0 & 0 & 0 \\
7 & 21 & 7 & 0 & 10 & 10 \\
  &    & 8 & 0 & 17 & 17 \\
  &    & 9 & 0 & 57 & 57 \\
  &    & 11 & 0 & 2 & 2 \\
7 & 22 & 8 & 0 & 11 & 11 \\
  &    & 9 & 0 & 90 & 90 \\
  &    & 11 & 0 & 9 & 9 \\
7 & 23 & $\lambda$ & 0 & 0 & 0 \\
7 & 24 & 10 & 0 & 45 & 45 \\
  &    & 11 & 2 & 3 & 3 \\
  &    & 12 & 2 & 55 & 55 \\
  &    & 13 & 0 & 1 & 1 \\
\hline
\end{tabular}
\vskip 5mm
\caption{Numbers of nonequivalent $\C^3_3(7,k,\lambda)$-cubes coming
from difference sets with multipliers.}\label{tab4}
\end{table}

The automorphism group $\Aut(C_7)$ is isomorphic
to~$C_6$. Automorphisms $\varphi\in \Aut(C_7)\setminus\{\text{id}\}$
keep only the neutral element of~$C_7$ fixed. Therefore, we can
restrict the search for $(7,k,\lambda)$ difference sets of
propriety~$3$ with a multiplier~$\varphi$ of order~$2$ or~$3$
to $k$-subsets of $G^2$ fixed by~$\varphi$, i.e.\ made up of
$\varphi$-orbits. The results of a computer enumeration up to
equivalence for $15\le k\le 24$ are summarized in Table~\ref{tab4}.
The columns labeled $C_7\rtimes C_2$ and $C_7\rtimes C_3$
contain the numbers of nonequivalent $(7,k,\lambda)$ cubes
constructed from difference sets of propriety~$3$ with multipliers
of order~$2$ and~$3$, respectively. The last column contains
the total number of nonequivalent cubes from the previous
two columns. For example, in the $(7,16,6)$-row there
is one cube in each of the columns $C_7\rtimes C_2$ and
$C_7\rtimes C_3$, and these two cubes are not equivalent.
In the $(7,19,9)$-row, there is only one cube constructed
from a difference set with a multiplier of order~$6$.

We now present constructions of difference sets of
propriety~$3$ from ordinary difference sets, producing
infinitely many examples. Let~$D$ be an ordinary
$(v,k,\lambda)$ difference set in a group~$G$. It is
known that the complement $\olsi{D}=G\setminus D$
is a difference set with parameters $(v,v-k,v-2k+\lambda)$.
From~\cite{RHB55}, it is also known that for every
$g\in G\setminus\{0\}$, there are exactly~$\lambda$
choices of $x,x'\in D$ with left difference $-x+x'=g$
(the definition involves right differences). We
will use this fact implicitly in the sequel. Moreover,
we will use

\begin{lemma}\label{lemmads}
For every $g\in G\setminus\{0\}$, there are exactly
$k-\lambda$ choices of $x\in D$, $x'\in \olsi{D}$
with $-x+x'=g$. The number of choices is the same
for $x\in \olsi{D}$, $x'\in D$.
\end{lemma}

Indeed, for fixed $g\ne 0$ and each $x\in D$, the equation $-x+x'=g$ determines $x'$ uniquely. Exactly~$\lambda$ of the resulting pairs have $x'\in D$, so the remaining $k-\lambda$ have $x'\in\olsi{D}$. The reverse case is analogous. Finally, we will use the following
notation for a Cartesian square without the diagonal:
\[ D^{2-}:=(D\times D) \setminus \{(x,x)\mid x\in D\}. \]

Difference sets with parameters $(4n-1,2n-1,n-1)$, $n\ge 2$
are known as \emph{Paley-type difference sets}.
The canonical examples are the nonzero squares in a finite
field of order $q\equiv 3 \pmod{4}$.
Many other series of Paley-type difference sets
are known, e.g.\ the twin prime power difference
sets~\cite[Theorem~VI.8.2]{BJL99}.

\begin{theorem}
Let $D\subseteq G$ be a Paley-type difference set
with parameters $(4n-1,2n-1,n-1)$. Then, $D^2\cup
(\olsi{D})^{2-} \subseteq G^2$
is a difference set of propriety~$3$ with parameters
\[(4n-1,\,(4n-1)(2n-1),\, (4n-1)(n-1)).\]
\end{theorem}

\begin{proof}
We use $(v,k,\lambda)$ and $(v,\olsi{k},\olsi{\lambda})$
for the parameters of~$D$ and~$\olsi{D}$. Clearly,
$K:=|D^2\cup (\olsi{D})^{2-}|=k^2+\olsi{k}(\olsi{k}-1)$.
Given~$g\in G\setminus\{0\}$, denote the number of
pairs $(x,y), (x',y')\in D^2\cup (\olsi{D})^{2-}$
satisfying conditions~1, 2, and~3 from Definition~\ref{ds3def}
by $\Lambda_1$, $\Lambda_2$, and $\Lambda_3$, respectively. We will express
these numbers in terms of the parameters. For the latter two,
direct counting shows that $\Lambda_2=\Lambda_3 = \lambda k+\olsi{\lambda}(\olsi{k}-2)$.
The former is expressed as
\[ \Lambda_1 = \lambda^2 + \olsi{\lambda}(\olsi{\lambda}-1)+
2(k-\lambda)(k-\lambda-1), \]
where Lemma~\ref{lemmads} is used for the last term.
Substituting $k=2n-1$, $\lambda=n-1$, $\olsi{k}=2n$, $\olsi{\lambda}=n$
gives $K=(4n-1)(2n-1)$, $\Lambda_1=\Lambda_2=\Lambda_3=(4n-1)(n-1)=:\Lambda$.
Therefore, $D^2\cup (\olsi{D})^{2-}$ is a difference set of
propriety~$3$ with parameters $(v,K,\Lambda)$.
\end{proof}

Another important class of difference sets are the \emph{Hadamard
difference sets} with parameters $(4n^2,2n^2-n,n^2-n)$, $n\ge 2$.
Many infinite series in this class are also known~\cite{DJ96}.

\begin{theorem}
Let $D\subseteq G$ be a Hadamard difference set
with parameters $(4n^2,2n^2-n,n^2-n)$. Then,
$D^2\cup (\olsi{D})^{2}$ is a
difference set of propriety~$3$ with parameters
$(4n^2,\, 2n^2(4n^2+1),\, 2n^2(2n^2+1))$. Moreover,
$D^{2-}\cup (\olsi{D})^{2-}$ and $(D\times \olsi{D})\cup
(\olsi{D}\times D)$ are difference sets of propriety~$3$
with parameters
$(4n^2,\, 2n^2(4n^2-1),\, 2n^2(2n^2-1))$.
\end{theorem}

\begin{proof}
We use the same notation as in the previous proof.
For the first set of pairs $D^2\cup (\olsi{D})^{2}$,
\begin{align*}
K &=k^2+\olsi{k}^2,\\
\Lambda_1&=\lambda^2+\olsi{\lambda}^2+2(k-\lambda)^2,\\
\Lambda_2&=\Lambda_3=\lambda k + \olsi{\lambda}\olsi{k}.
\end{align*}
Similarly, for the second set of pairs $D^{2-}\cup (\olsi{D})^{2-}$,
\begin{align*}
K &= k(k-1)+\olsi{k}(\olsi{k}-1),\\
\Lambda_1 &= \lambda(\lambda-1)+\olsi{\lambda}(\olsi{\lambda}-1)+ 2(k-\lambda)(k-\lambda-1),\\
\Lambda_2 &=\Lambda_3 = \lambda (k-2) + \olsi{\lambda}(\olsi{k}-2),
\end{align*}
and for the third set of pairs $(D\times \olsi{D})\cup (\olsi{D}\times D)$,
\begin{align*}
K &= 2k \olsi{k},\\
\Lambda_1 &= 2\lambda \olsi{\lambda} + 2(k-\lambda)^2,\\
\Lambda_2 &=\Lambda_3 = \lambda \olsi{k} + \olsi{\lambda}k.
\end{align*}
Substituting $k=2n^2-n$, $\lambda=n^2-n$, $\olsi{k}=2n^2+n$, $\olsi{\lambda}=n^2+n$
shows that $\Lambda_1=\Lambda_2=\Lambda_3$ and gives parameters for difference sets
of propriety~$3$ as stated in the theorem.
\end{proof}

\section{Infinite families of three-dimensional symmetric designs}\label{sec4}

We first describe constructions producing infinite families
of $\C_3^3(v,k,\lambda)$-cubes inspired by known constructions
of higher-dimensional Hadamard matrices. Let~$A$ be a Hadamard
matrix of order~$v$, dimension and propriety~$n$, and~$L$ a
Latin square of order~$v$. By~\cite[Theorem~6.1.8]{YNX10}, the $(n+1)$-dimensional array $B$ with
$$B(x_1,\ldots,x_n,x_{n+1}) := A\big(x_1,\ldots,x_{n-1},L(x_n,x_{n+1})\big), \quad \mbox{ for } x_1,\dots,x_{n+1}\in \{1,\dots,v\}$$
is a Hadamard matrix of order~$v$ and propriety $n+1$.
A  similar construction applies to symmetric designs of dimension
and propriety~$n$, illustrated here for $n=2$.

\begin{theorem}\label{tmdimpp}
Let $A$ be an incidence matrix of a symmetric $(v,k,\lambda)$ design,
and let $L$ be a Latin square of order~$v$. Then, the array $B$ defined as follows is a $\C_3^3(v,vk,v\lambda)$-cube:
$$B(x,y,z)=A\big(x,L(y,z)\big), \quad \mbox{ for } x,y,z\in \{1,\dots,v\}.$$
\end{theorem}

\begin{proof}
We show that~$B$ satisfies one of the conditions in ~\eqref{maincond}; the proofs of the other conditions are very similar. Let us fix  $y_1,y_2\in \{1,\ldots,v\}$.  If $y_1=y_2$, since $L$ is a Latin square, we have $\{L(y_1,z)\ |\ 1\leq z\leq v\}=\{1,\dots,v\}$, and so $\sum_{z=1}^v  A\big(x,L(y_1,z)\big)$ is the number of ones in row~$x$ of~$A$, for $x\in \{1,\dots, v\}$.
If $y_1\neq y_2$, then $c_1:=L(y_1,z)\neq L(y_2,z)=:c_2$, and so $\sum_{x=1}^v A(x,c_1)A(x,c_2)$ is the scalar product of columns $c_1$ and $c_2$ of $A$, for $z\in \{1,\dots, v\}$. The following completes the proof:
\begin{equation*}
\sum_{z=1}^v\sum_{x=1}^v  A\big(x,L(y_1,z)\big)A\big(x,L(y_2,z)\big)=
\begin{cases}
    \displaystyle\sum_{x=1}^v\sum_{z=1}^v  A\big(x,L(y_1,z)\big)=\sum_{x=1}^v k = vk, &\mbox{ if }y_1=y_2, \\[5mm]
    \displaystyle\sum_{z=1}^v \lambda  =v\lambda, &\mbox{ if }y_1\neq y_2. %\qedhere
\end{cases}
\end{equation*}
\end{proof}

\begin{remark}\textup{
Permuting the rows of~$A$ and the rows and columns of~$L$ in the previous
theorem clearly gives isotopic cubes~$B$. However, permuting the
columns of~$A$ or the symbols of~$L$ generally gives nonequivalent
cubes, as well as using conjugates of~$L$ with exchanged symbols and rows
or columns. For example, by a computation in PAG~\cite{PAG}, the~$7!$
column-permutations of an incidence matrix of the Fano plane and representatives
of the $564$ isotopy classes of Latin squares of order~$7$ give rise to
$7417$ nonequivalent cubes in $C^3_3(7,21,7)$. One of these cubes
is a $C^3_2(7,3,1)$-cube of propriety~$2$; it can be obtained from a circulant
incidence matrix of the Fano plane and the cyclic Latin square of
order~$7$:
$$A=\left[
\begin{array}{lllllll}
 1 & 1 & 0 & 1 & 0 & 0 & 0 \\
 0 & 1 & 1 & 0 & 1 & 0 & 0 \\
 0 & 0 & 1 & 1 & 0 & 1 & 0 \\
 0 & 0 & 0 & 1 & 1 & 0 & 1 \\
 1 & 0 & 0 & 0 & 1 & 1 & 0 \\
 0 & 1 & 0 & 0 & 0 & 1 & 1 \\
 1 & 0 & 1 & 0 & 0 & 0 & 1
\end{array}
\right] \quad \mbox{ and } \quad
L=\left[
\begin{array}{lllllll}
 1 & 2 & 3 & 4 & 5 & 6 & 7 \\
 2 & 3 & 4 & 5 & 6 & 7 & 1 \\
 3 & 4 & 5 & 6 & 7 & 1 & 2 \\
 4 & 5 & 6 & 7 & 1 & 2 & 3 \\
 5 & 6 & 7 & 1 & 2 & 3 & 4 \\
 6 & 7 & 1 & 2 & 3 & 4 & 5 \\
 7 & 1 & 2 & 3 & 4 & 5 & 6
\end{array}
\right].$$
In general, a cube of propriety~$2$ is obtained if~$L$ is the multiplication
table of a group~$G$, and the rows of~$A$ are incidence vectors of $(v,k,\lambda)$
difference sets in~$G$. This is the $n=3$ case of the ``group cube'' construction
\cite[Theorem~4.1]{KPT25}.
}\end{remark}

\begin{remark}\textup{
In~\cite[Problem~3]{KZXS25}, it was asked whether $\C^n_2(25,9,3)$-cubes of
propriety~$2$ exist for $n\ge 3$. In the three-dimensional case, such cubes
would also be $\C^3_3(25,25\cdot 9,25\cdot 3)$-cubes. Choose in Theorem~\ref{tmdimpp} the multiplication table~$L$ of a group~$G$ of order~$25$. The resulting cubes of propriety~$3$ are not of
propriety~$2$: otherwise,  the supports of the rows of~$A$ would be $(25,9,3)$ difference sets in~$G$, but no such difference sets exist in groups of order~$25$.
}\end{remark}

The next lemma will be used in the proofs of the upcoming constructions.
For $\alpha\in \mathbb{Z}, \beta \in \mathbb{N}\cup \{0\}$, a $v\times v\times v$ array $C$  is called \emph{$(\alpha,\beta)$-nice} if each entry of $C$ is $+1$ or $-1$, the sum of the entries of each layer of~$C$ is~$\alpha$, and the scalar product of every two distinct parallel layers is~$\beta$.

\begin{lemma}\label{lem:symcube}
If $C$ is an $(\alpha,\beta)$-nice $v\times v\times v$ array, then $(J- C)/2$ is a $\C_3^3(v,\frac{v^2- \alpha}{2},\frac{v^2-2\alpha+\beta}{4})$-cube.
\end{lemma}

\begin{proof}
    If we replace each $+1$ in $C$  by a 0, and each $-1$ in $C$ by a 1, we obtain $B:=(J- C)/2$. Let $L$ be a layer in $C$, and let $\oplus$ and $\ominus$ be the number of occurrences of $+1$ and $-1$ in $L$. Since $\oplus +\ominus=v^2$ and $\oplus -\ominus=\alpha$, we conclude that the number of occurrences of $1$ in the layer of $B$ that corresponds to $L$ is  $\ominus=(v^2-\alpha)/2$.  Now, let $L'$ be another layer of $C$ that is parallel to $L$. Let $\boxplus$ and $\boxminus$ be the number of positions where both $L$ and $L'$ are $+1$ and $-1$, respectively.  Since the scalar product of $L$ and $L'$   is $\beta$, we have $\beta = \boxplus +\boxminus - (v^2-\boxplus-\boxminus)$. But  $\boxplus-\boxminus=\oplus -\ominus=\alpha$,  and so the scalar product of the two layers  of $B$ that correspond to $L$  and $L'$ is $\boxminus=(v^2-2\alpha+\beta)/4$.
\end{proof}

\subsection{Product constructions}

Another construction of $n$-dimensional Hadamard matrices is
the so-called \emph{product construction}, proved independently by
Yang~\cite{YXY86} and de~Launey~\cite{WdL87} (see
also~\cite[Section~5.1.2]{KH07} and~\cite[Theorem~6.1.19]{YNX10}).
Given a Hadamard matrix~$H$ of order~$v$, the
  $n$-dimensional  matrix $A$ with
$$A(x_1,\ldots,x_n) := \prod_{1\le i < j \le n} H(x_i,x_j), \quad \mbox{ for } x_1,\dots,x_{n}\in \{1,\dots,v\}$$
is a Hadamard matrix of order~$v$ and propriety~$2$.
An analogous construction for symmetric designs gives cubes
of propriety~$3$ and relies on the following lemma, where a matrix  whose diagonal entries
are all $i$ is called \emph{$i$-diagonal}.

\begin{lemma}
Every nontrivial symmetric design possesses both a $0$-diagonal incidence matrix,
and a $1$-diagonal incidence matrix.
\end{lemma}

\begin{proof}
Let us assume that the rows and columns of an incidence matrix~$A$ of a nontrivial
symmetric design~$\mathcal{A}$ are indexed by the points $x_1,\ldots,x_v$ and the
blocks $y_1,\ldots,y_v$. Let~$\Gamma$ be the bipartite graph whose
vertices are the points and blocks of~$\mathcal{A}$ such that~$x_i$ and~$y_j$ are
adjacent if $x_i\in y_j$. By Hall's marriage theorem, since~$\Gamma$ is $k$-regular,
it has a perfect matching in which $x_i$ is matched with $y_{\pi(i)}$ for some
permutation $\pi\in S_v$.  The matrix $\pi(A)$ whose column~$i$ is column~$\pi(i)$
of~$A$, for $i\in \{1,\dots,v\}$, is a $1$-diagonal incidence matrix of~$\mathcal{A}$.
Repeating the same argument for the complement of~$\Gamma$ gives a $0$-diagonal
incidence matrix of~$\mathcal{A}$.
\end{proof}

\begin{theorem}\label{thmprodconst}
Let $A$ be an $i$-diagonal incidence matrix of a symmetric $(v,k,\lambda)$ design, where $i\in \{0,1\}$. Then, the array $B$ defined by
$$B(x,y,z)=A(x,y)+A(y,z) +A(x,z)\,\bmod{\,2}, \quad \mbox{ for } x,y,z\in \{1,\dots,v\}$$
is a $\C_3^3(v,K_i,\Lambda_i)$-cube, where
\begin{alignat*}{3}
K_0 &=k (3 v + 4\lambda -6 k ), & \quad \Lambda_0 &=k ( v - 4 \lambda-2 ) + 2 \lambda (2 \lambda+1),\\[1mm]
K_1 &=k (3 v + 4\lambda -2 k) - 4 v\lambda, & \quad \Lambda_1 &= k (v - 4 \lambda+2) + 2 \lambda (2 \lambda-1).
\end{alignat*}
\end{theorem}

\begin{proof}
Let  $C$ be the $v\times v\times v$ array over $\{-1,1\}$ with
\[
C(x,y,z)=\big(1-2A(x,y)\big)\big(1-2A(y,z)\big)\big(1-2A(x,z)\big), \quad \mbox{ for }x,y,z\in \{1,\dots,v\}.
\]
We remark that $B=(J-C)/2$. The sum of the entries of each layer obtained by
fixing the first coordinate of $C$ is as follows:
\begin{align*}
\alpha_i&:=\sum_{y=1}^v \sum_{z=1}^v \big(1-2A(x,y)\big)\big(1-2A(x,z)\big)\big(1-2A(y,z)\big) \\
 &=\sum_{y=1}^v  \big(1-2A(x,y)\big) \sum_{z=1}^v \big(1-2A(x,z)-2A(y,z)+4A(x,z)A(y,z)\big) \\
 &=\sum_{y=1}^v \big(1-2A(x,y)\big) \big(v-4k+4[x=y]k+4[x\neq y]\lambda\big) \\
 &=\sum_{y=1}^v \Big( \big(v-4k\big)\big(1-2A(x,y)\big) + 4\big(1-2A(x,y)\big)\big([x=y]k+[x\neq y]\lambda\big)\Big) \\
 &= (v-4k)(v-2k)+4k\big(1-2A(x,x)\big)+4\lambda\big((v-1)-2\big(k-A(x,x)\big)\big) \\
 &=(v-4k)(v-2k)+4k+4\lambda(v-2k-1)+ 8(\lambda-k)i.
\end{align*}

The scalar product of two distinct layers obtained by fixing the first coordinate to two distinct values is
\begin{align*}
\beta&:=\sum_{y=1}^v \sum_{z=1}^v  \big(1-2A(x_1,y)\big)\big(1-2A(x_2,y)\big)\big(1-2A(x_1,z)\big)\big(1-2A(x_2,z)\big)\big(1-2A(y,z)\big)^2 \\
&= \sum_{y=1}^v \big(1-2A(x_1,y)\big)\big(1-2A(x_2,y)\big)\sum_{z=1}^v \big(1-2A(x_1,z)\big)\big(1-2A(x_2,z)\big) \\
&=\Bigg(\sum_{y=1}^v \big(1-2A(x_1,y)\big)\big(1-2A(x_2,y)\big)\Bigg)^2 \\
&=\Bigg(\sum_{y=1}^v \big(1-2A(x_1,y)-2A(x_2,y)+4A(x_1,y)A(x_2,y)\big)\Bigg)^2 =(v-4k+4\lambda)^2.
\end{align*}
The argument is very similar if the second or third coordinates are fixed.
Thus, $C$ is $(\alpha_i,\beta)$-nice and so, by Lemma~\ref{lem:symcube}, $B$ is a $\C_3^3(v,\frac{v^2- \alpha_i}{2},\frac{v^2-2\alpha_i+\beta}{4})$-cube. Using $\lambda(v-1)=k(k-1)$
for simplification, the formulae for $K_i$ and $\Lambda_i$ from the statement of the
theorem are obtained.
\end{proof}

\begin{remark}\textup{
Surprisingly, the construction of Theorem~\ref{thmprodconst}
may produce nonequivalent cubes even from isomorphic incidence
matrices. For example, consider the following incidence matrices of
the Fano plane:
\[
A_1 =
\begin{bmatrix}
0 & 0 & 0 & 1 & 1 & 1 & 0 \\
1 & 0 & 1 & 0 & 0 & 1 & 0 \\
1 & 0 & 0 & 0 & 1 & 0 & 1 \\
0 & 1 & 1 & 0 & 1 & 0 & 0 \\
0 & 1 & 0 & 0 & 0 & 1 & 1 \\
0 & 0 & 1 & 1 & 0 & 0 & 1 \\
1 & 1 & 0 & 1 & 0 & 0 & 0
\end{bmatrix}, \quad
A_2 =
\begin{bmatrix}
0 & 0 & 0 & 1 & 1 & 0 & 1 \\
1 & 0 & 1 & 1 & 0 & 0 & 0 \\
1 & 0 & 0 & 0 & 1 & 1 & 0 \\
0 & 1 & 1 & 0 & 1 & 0 & 0 \\
0 & 0 & 1 & 0 & 0 & 1 & 1 \\
1 & 1 & 0 & 0 & 0 & 0 & 1 \\
0 & 1 & 0 & 1 & 0 & 1 & 0
\end{bmatrix}.
\]
The full autotopy groups of the corresponding cubes are
$\Atop(B_1)\cong C_7\rtimes C_3$ and $\Atop(B_2)\cong C_3$, therefore
$B_1$ and $B_2$ are not equivalent. In~\cite{dLS08}, it was shown that applying the product construction to Hadamard matrices equivalent to the Sylvester matrix of order~$2^n$
results in exponentially many nonequivalent Hadamard cubes, thereby identifying an error  in~\cite{KM02}.
}\end{remark}

A \emph{doubly regular tournament} of order $4t-1$ is a symmetric $(4t-1,2t-1,t-1)$ design whose incidence matrix $A$ satisfies $A+A^\top=J-I$. Doubly regular tournaments are equivalent
to skew Hadamard matrices of order~$4t$~\cite{RB72}, and are conjectured to exist
for all $t\ge 1$~\cite[Conjecture~1.21]{CK07}. The next construction is different from Theorem~\ref{thmprodconst}: its third term is $A_i(z,x)$ rather than $A_i(x,z)$, and the tournament identity for $A_i^2$ is essential.

\begin{theorem}\label{thmprodconst2}
Let $A_0$ be the incidence matrix of a doubly regular tournament of order $4t-1$, and $A_1=A_0+I$, where $I$ is the identity matrix. Then the array $B_i$, $i\in\{0,1\}$, defined by
$$B_i(x,y,z)=A_i(x,y)+A_i(y,z) +A_i(z,x)\,\bmod{\,2}, \quad \mbox{ for } x,y,z\in \{1,\dots,4t-1\}$$
is a $\C_3^3(4t-1,K_i,\Lambda_i)$-cube, where
\begin{alignat*}{3}
K_0 &=(2t-1)(4t+3), & \quad \Lambda_0 &=(2t-1)(2t+3),\\[1mm]
K_1 &=2(4t^2-5t+2), & \quad \Lambda_1 &= 4(t-1)^2.
\end{alignat*}
\end{theorem}

\begin{proof}
Let $v=4t-1$, $k=2t-1$, $\lambda=t-1$, and $k_i=k+i$, $\lambda_i=\lambda+i$, where $i\in\{0,1\}$. Let~$C_i$ be the $v\times v\times v$ array over $\{-1,1\}$ with
\[
C_i(x,y,z)=\big(1-2A_i(x,y)\big)\big(1-2A_i(y,z)\big)\big(1-2A_i(z,x)\big), \quad \mbox{ for }x,y,z\in \{1,\dots,v\}.
\]
We remark that $B_i=(J-C_i)/2$, $A_iA_i^\top=(k-\lambda)I+\lambda_i J$, and
$A_i^2=t(J-I)-(-1)^{i}A_i$.  The sum of the entries of each layer obtained
by fixing the first coordinate of $C_i$ is as follows:
\begin{align*}
\alpha_i&:=\sum_{y=1}^v \sum_{z=1}^v \big(1-2A_i(x,y)\big)\big(1-2A_i(y,z)\big)\big(1-2A_i(z,x)\big) \\
 &=\sum_{y=1}^v  \big(1-2A_i(x,y)\big) \sum_{z=1}^v \big(1-2A_i(z,x)-2A_i(y,z)+4A_i(y,z)A_i(z,x)\big) \\
  &=\sum_{y=1}^v  \big(1-2A_i(x,y)\big) \big(v-4k_i+4A_i^2(y,x)\big) \\
 &=\sum_{y=1}^v  \big(1-2A_i(x,y)\big) \big(v-4k_i+4t[x\neq y]-(-1)^i4A_i(y,x)\big) \\
&=\sum_{y=1}^v \Big(\big(1-2A_i(x,y)\big)\big(v-4k_i+4t[x\neq y]\big)-(-1)^i  \big(4A_i(y,x) -8A_i(x,y)A_i(y,x) \big) \Big)\\
&=\sum_{y=1}^v \Big( v-4k_i+4t[x\neq y]-2A_i(x,y)\big(v-4k_i+4t[x\neq y]\big)-(-1)^i  \big(4A_i(y,x) -8[x=y]i \big) \Big) \\
&=v(v-4k_i)+4t(v-1)-2k_i(v-4k_i)-8tk-(-1)^i(4k_i-8i) = (-1)^i(7-12t).\\
\end{align*}
The scalar product of two distinct layers obtained by fixing the first coordinate to two distinct values is
\begin{align*}
\beta
&:= \sum_{y=1}^v \sum_{z=1}^v
   \big(1-2A_i(x_1,y)\big)\big(1-2A_i(x_2,y)\big)
   \big(1-2A_i(z,x_1)\big)\big(1-2A_i(z,x_2)\big)\big(1-2A_i(y,z)\big)^2 \\
&= \sum_{y=1}^v
      \big(1-2A_i(x_1,y)\big)\big(1-2A_i(x_2,y)\big)
   \sum_{z=1}^v
      \big(1-2A_i(z,x_1)\big)\big(1-2A_i(z,x_2)\big) \\
&= \sum_{y=1}^v
      \big(1-2A_i(x_1,y)-2A_i(x_2,y)+4A_i(x_1,y)A_i(x_2,y)\big)  \\
&\quad\;\;\times
   \sum_{z=1}^v
      \big(1-2A_i(z,x_1)-2A_i(z,x_2)+4A_i(z,x_1)A_i(z,x_2)\big) = (\,v - 4k_i + 4\lambda_i\,)^2=1.
\end{align*}
The argument is very similar if the second or third coordinates are fixed.
Thus, $C_i$ is $(\alpha_i,\beta)$-nice and so, by Lemma~\ref{lem:symcube}, $B_i$ is a $\C_3^3(v,\frac{v^2- \alpha_i}{2},\frac{v^2-2\alpha_i+\beta}{4})$-cube.
\end{proof}

\begin{remark}\textup{
Equivalence of doubly regular tournaments is defined by applying the same permutation
to the rows and columns of~$A$, in order to preserve the property $A+A^\top=J-I$. Such
incidence matrices~$A_i$ will clearly produce equivalent cubes $B_i$ from
Theorem~\ref{thmprodconst2}. However, applying the construction to doubly regular
tournaments which are not equivalent generally gives nonequivalent cubes.
For example, in~\cite{HKMTR20} it was shown that there are precisely $7227$ skew
Hadamard matrices of order $32$ up to SH-equivalence. These matrices can be downloaded
from~\cite{TRweb}. A computation in PAG~\cite{PAG} showed that $7061$ nonequivalent cubes
in each of the sets $\C^3_3(31,525,285)$ and $\C^3_3(31,436, 196)$ are obtained by applying Theorem~\ref{thmprodconst2} to the corresponding doubly regular tournaments of order~$31$.
Many more nonequivalent cubes in $\C^3_3(31,465, 225)$, $\C^3_3(31,464, 224)$,
$\C^3_3(31,497, 257)$, and $\C^3_3(31,496, 256)$ can be constructed
by independently permuting the rows and columns of $A$ and $A+I$, putting
them in $0$-diagonal or $1$-diagonal form, and applying Theorem~\ref{thmprodconst}.
}\end{remark}

\subsection{Constructions using Hadamard matrices and layer-rainbow cubes}

Our next two constructions of $\C^3_3$-cubes use Hadamard
matrices of square order and another kind of Latin cubes.
A $v\times v\times v$ array on $v^2$ symbols is a
\emph{layer-rainbow cube} if each layer contains each symbol.
A more general situation is dealt with in~\cite{BS}. Layer-rainbow
cubes of order $v$ are equivalent to one-factorizations of complete 3-partite 3-uniform hypergraphs \cite{BSid23,BProc24}.  Factorizations of complete multipartite hypergraphs are known; see~\cite{Baran79}. Latin cubes and hypercubes were studied early by Kishen~\cite{KK49}. The following elementary lemma provides a simple construction of layer-rainbow cubes of all orders~$v$.

\begin{lemma}\label{latinprod}
If $L_1$ and $L_2$ are two Latin squares of order~$v$
with entries in $\{1,\ldots,v\}$, then
$L(x,y,z)=v(L_1(x,y)-1)+L_2(y,z)$
defines a layer-rainbow  cube of order~$v$
with entries in $\{1,\ldots,v^2\}$.
\end{lemma}

A Hadamard matrix is \emph{normalized} if its first row consists entirely of ones, and \emph{regular} if all rows and columns have the same sum. In that case, the order of
the Hadamard matrix is~$v^2$, and all row- and column sums are~$v$, or all sums are~$-v$.
Multiplying by~$-1$ if necessary, we can assume that all sums are~$v$.

\begin{theorem}\label{tmhadsqconstr}
%Let $H$ be a Hadamard matrix of order $v^2$.
Assume a Hadamard matrix $H$ of order $v^2$ exists.
\begin{enumerate}
\item[\textup{(i)}]  There exists
a $\C^3_3(v^3,v^4(v^2-1)/2,v^4(v^2-2)/4)$-cube.
\item[\textup{(ii)}] If $H$ is regular, there exists a
$\C^3_3(v^3,v^3(v^3-1)/2,v^3(v^3-2)/4)$-cube.
\end{enumerate}
\end{theorem}

\begin{proof}
Let $L$ be a layer-rainbow cube of order $v$ with entries in $\{1,\ldots,v^2\}$, and define   the sequence $M_1,\dots,M_{v^2}$  of $v^2\times v^2 \times v^2$ arrays as follows:
$$M_{\ell}(x,y,z)=H(\ell,x)H(\ell,y)H(\ell,z), \quad \mbox{ for }\ell,x,y,z\in\{1,\ldots,v^2\}.$$
If $H$ is normalized, $M_1=J$.
Let $C$ be a  $v^3\times v^3\times v^3$  array over $\{1,-1\}$, formed by replacing each occurrence of $\ell$ in $L$ with the array $M_{\ell}$, for $\ell \in \{1,\ldots,v^2\}$.  We claim that  $C$ is $(\alpha,0)$-nice, where $\alpha=v^4$ if $H$ is normalized, and $\alpha=v^3$ if $H$ is regular.
Consequently, by Lemma~\ref{lem:symcube}, $(J- C)/2$ is a $\C_3^3(v^3,\frac{v^6- \alpha}{2},\frac{v^6-2\alpha}{4})$-cube.

Let $\alpha$ be the  sum of the entries  of each layer of $C$ obtained by fixing the first coordinate. Note that
\begin{align*}
\gamma_{\ell}(x):=\sum_{y=1}^{v^2}\sum_{z=1}^{v^2}M_{\ell}(x,y,z)= H(\ell,x)\sum_{y=1}^{v^2}H(\ell,y)\sum_{z=1}^{v^2}H(\ell,z) =
    \begin{cases}
        v^4 [\ell=1], &\text{if $H$ is normalized},\\
         v^2 H(\ell,x), &\text{if $H$ is regular}.
    \end{cases}
\end{align*}
But $\alpha=\sum_{\ell=1}^{v^2} \gamma_{\ell}(x)$, and so $\alpha=v^4$ if $H$ is normalized, and $\alpha=v^3$ if $H$ is regular.
A similar argument applies to layers with the second or third coordinate fixed.

To complete the proof, we  calculate $\beta:=\sum_{y=1}^{v^3}\sum_{z=1}^{v^3}C(x,y,z)C(\dot x,y,z)$ for distinct $x,\dot x\in \{1,\ldots,v^3\}$. The argument applies to all Hadamard matrices, irrespective of whether they are normalized or regular, and proceeds analogously when the second or third coordinate is held fixed. Let $w$ be any one of the coordinates $x$, $\dot x$, $y$, or $z$. We can write $w$ as $w=v^2(w_{1}-1)+w_{2}$ for some $w_{1}\in\{1,\ldots,v\}$, $w_2\in\{1,\ldots,v^2\}$. Then, we have
\begin{align*}
\delta(y_1,y_2,z_1,z_2)&:=M_{L(x_{1},y_1,z_1)}\big(x_2,y_2,z_2\big) M_{L(\dot x_{1},y_1,z_1)}\big(\dot x_2,y_2,z_2\big)\\
&=\prod_{u\in\{x_2,y_2,z_2\}} H\big( L(x_{1},y_1,z_1),u  \big) \prod_{u\in\{\dot x_2,y_2,z_2\}} H\big( L(\dot x_{1},y_1,z_1),u  \big).
\end{align*}
If $x_1=\dot x_1$, then $x_2\neq \dot x_2$, and since $L$ is layer-rainbow, we have
\begin{align*}
\sum_{y_1=1}^v \sum_{z_1=1}^v \delta(y_1,y_2,z_1,z_2)
=  \sum_{\ell=1}^{v^2} H(\ell,x_{2})H(\ell,\dot x_{2})=0.
\end{align*}
If $x_1\neq \dot x_1$, then $\ell:=L(x_{1},y_1,z_1)\neq L(\dot x_{1},y_1,z_1)=:\dot \ell $, and so
\begin{align*}
&\sum_{y_2=1}^{v^2}\sum_{z_2=1}^{v^2}\delta(y_1,y_2,z_1,z_2)
=H(\ell,x_{2}) H(\dot \ell ,\dot x_{2}) \left (\sum_{y_2=1}^{v^2} H(\ell,y_2)H(\dot \ell ,y_2)\right)
  \left (\sum_{z_2=1}^{v^2} H( \ell ,z_2)H(\dot \ell  , z_2)\right) =0.
\end{align*}
But  $\beta=\sum_{\genfrac{}{}{0pt}{}{y_i,z_i\in \{1,\dots,v^i\}}{i\in \{1,2\}}}\delta(y_1,y_2,z_1,z_2)$, and so in both cases we have $\beta=0$.
\end{proof}

\begin{comment}
\begin{remark}\textup{
Row-sums of the array~$C$ are of the following form:
    \begin{align*}
        & \sum_{i\in \Gamma}\sum_{z=1}^{v^2}M_i(a,b,z) =\sum_{i\in \Gamma}H(i,a)H(i,b)\sum_{z=1}^{v^2}H(i,z) =\\[2mm]
    &=\begin{cases}
        v^2\displaystyle\sum_{i\in \Gamma}H(i,a)H(i,b)[i=1] &\text{if $H$ is normalized},\\[1mm]
        v\displaystyle\sum_{i\in \Gamma}H(i,a)H(i,b) &\text{if $H$ is regular}, \\
    \end{cases}\\[2mm]
    &=\begin{cases}
        v^2H(1,a)H(1,b)[i\in \Gamma] &\text{if $H$ is normalized},\\[1mm]
        v\displaystyle\sum_{i\in \Gamma}H(i,a)H(i,b) &\text{if $H$ is regular},
    \end{cases}
    \end{align*}
where $\Gamma=\Gamma_{c,d}=\{L(c,d,i) \mid i\in\{1,\ldots,v\}\}$ for $c,d\in\{1,\ldots,v\}$.
For (i), a row-sum of $C$ is $0$ if $i\in \Gamma,a=b$ and $\pm v^2$ otherwise.
For (ii), the row-sum corresponding to $a=b$ is $v^2$, and it is not difficult to see
that for some distinct $a,b$ and some $c,d$, the row-sum cannot be $v^2$.
Therefore, the resulting three-dimensional symmetric designs are not of propriety $2$.
}\end{remark}
\end{comment}

\begin{remark}\textup{
Generally, using nonequivalent Hadamard matrices~$H$ and/or
nonequivalent layer-rainbow Latin cubes~$L$ in
Theorem~\ref{tmhadsqconstr} produces nonequivalent
$\C_3^3$-cubes. For example, there are five nonequivalent
Hadamard matrices of order~$16$~\cite{MH61, NJAS}. Using
Lemma~\ref{latinprod}, three nonequivalent layer-rainbow
Latin cubes of order~$4$ are obtained as products of
Latin squares from the two main classes of order~$4$.
By a computation in PAG~\cite{PAG}, the $15$ cubes in
$\C_3^3(64, 1920, 896)$ resulting from Theorem~\ref{tmhadsqconstr}~(i)
are pairwise nonequivalent. Three of the five equivalence
classes of Hadamard matrices of order~$16$ contain regular
Hadamard matrices. The $9$ cubes in $\C_3^3(64,2016,992)$ resulting
from Theorem~\ref{tmhadsqconstr}~(ii) are also pairwise nonequivalent.
}\end{remark}

We can construct $\C^3_3$-cubes exhibiting additional regularities
by using special layer-rainbow cubes~$L$ in Theorem~\ref{tmhadsqconstr}.
A three-dimensional symmetric design $A$ of order $v^3$ is
said to be \emph{block-symmetric} if there is a partition of $A$ into
$v^2\times v^2\times v^2$ arrays $A_{ijk}$, $i,j,k\in \{1,\dots,v\}$ such that
\begin{enumerate}
    \item $A_{ij\ell}=A_{j\ell i}=A_{\ell ij}$ for distinct $i,j,\ell \in  \{1,\dots,v\}$,
    \item  $A_{iij}=A_{jj i}, A_{iji}=A_{jij}, A_{ij j}=A_{jii}$ for  $i,j\in \{1,\dots,v\}$.
\end{enumerate}

\begin{corollary}\label{tmhadsqconstrSYM}
Let $H$ be a Hadamard matrix of order $v^2$ and $v\equiv 0,2 \pmod{3}$, $v\neq 3$.
\begin{enumerate}
\item[\textup{(i)}]  There exists
a block-symmetric $\C^3_3(v^3,v^4(v^2-1)/2,v^4(v^2-2)/4)$-cube.
\item[\textup{(ii)}] If $H$ is regular, there exists a block-symmetric $\C^3_3(v^3,v^3(v^3-1)/2,v^3(v^3-2)/4)$-cube.
\end{enumerate}
\end{corollary}

\begin{proof}
By \cite[Theorem~1.1]{BSid23}, there exists a \emph{symmetric} layer-rainbow cube~$L$ of order $v$, with the following properties:
\begin{enumerate}
\item $L(x,y,z)=L(y,z,x)=L(z,x,y)$ for distinct $x,y,z \in  \{1,\dots,v\}$,
\item $L(x,x,y)=L(y,y,x), L(x,y,x)=L(y,x,y), L(x,y,y)=L(y,x,x)$ for  $x,y\in \{1,\dots,v\}$.
\end{enumerate}
Using~$L$ in the proof of Theorem~\ref{tmhadsqconstr} gives block-symmetric $\C^3_3$-cubes.
\end{proof}

\begin{comment}
    For example, the following symmetric layer-rainbow cube of order~$6$ can be used together
with Hadamard matrices of order~$36$ to construct block-symmetric $\C^3_3$-cubes of order~$216$.
\vskip 3mm

\begin{center}
\tabcolsep=0.031cm
\begin{tabular}{|l|l|l|l|l|l|}
\hline
\begin{tabular}{ c c c c c c}
 j & D & B & F & K & A \\
 C & M & U & h & P & W \\
 L & V & H & d & T & Y  \\
 O & R & a & I & i & c \\
 N & Z & b & Q & J & X \\
 E & f & S & g & e & G
\end{tabular}
&
\begin{tabular}{ c c c c c c}
 M & C & V & R & Z & f \\
 D & j & A & J & B & F \\
 U & I & O & X & g & Q  \\
 h & N & e & G & S & T \\
 P & H & c & Y & E & a \\
 W & L & i & b & d & K
\end{tabular}
&
\begin{tabular}{ c c c c c c}
 H & U & L & a & b & S \\
 V & O & I & e & c & i \\
 B & A & j & M & C & D  \\
 d & X & E & K & f & Z \\
 T & g & F & W & G & h \\
 Y & Q & J & P & R & N
\end{tabular}
&
\begin{tabular}{ c c c c c c}
 I & h & d & O & Q & g \\
 R & G & X & N & Y & b \\
 a & e & K & E & W & P  \\
 F & J & M & j & D & C \\
 i & S & f & A & L & U \\
 c & T & Z & B & V & H
\end{tabular}
&
\begin{tabular}{ c c c c c c}
 J & P & T & i & N & e \\
 Z & E & g & S & H & d \\
 b & c & G & f & F & R  \\
 Q & Y & W & L & A & V \\
 K & B & C & D & j & O \\
 X & a & h & U & I & M
\end{tabular}
&
\begin{tabular}{ c c c c c c}
 G & W & Y & c & X & E \\
 f & K & Q & T & a & L \\
 S & i & N & Z & h & J  \\
 g & b & P & H & U & B \\
 e & d & R & V & M & I \\
 A & F & D & C & O & j
\end{tabular}
\\ \hline
\end{tabular}
\end{center}
\vskip 2mm
\end{comment}

The construction of Theorem~\ref{tmhadsqconstr}~(i) also gives rise to
$\C^3_3$-cubes with large subcubes. Wallis~\cite{Wal86} showed that if~$H$
is any $v\times v$ matrix over $\{1,-1\}$, and if~$w$ is the smallest
number no less than~$v$ which is the order of a Hadamard matrix, then~$H$
can be embedded in a Hadamard matrix of order~$w^2$. In particular,
if~$H$ is a Hadamard matrix of order~$v$, then there is a Hadamard matrix
of order~$v^2$ containing $H$ as a submatrix.

\begin{theorem} \label{embsymcub}
Let $H$ be a Hadamard matrix of order $v^2$.
There exists a $\C^3_3(v^6,v^8(v^4-1)/2,v^8(v^4-2)/4)$-cube containing
a $\C^3_3(v^3,v^4(v^2-1)/2,v^4(v^2-2)/4)$-cube as a subcube.
\end{theorem}

\begin{proof}
By Wallis' theorem, there exists a Hadamard matrix~$K$ of order~$v^4$ containing~$H$.
After permuting rows and columns so that a row and a column of the embedded~$H$ come first, normalize~$K$; the embedded copy of~$H$ is then normalized as well. Let $L_1$ be a layer-rainbow
cube of order $v$ with entries in $\{1,\ldots,v^2\}$.  Applying  Theorem~\ref{tmhadsqconstr}~(i)
with~$H$ and~$L_1$, we obtain a~$\C^3_3(v^3,v^4(v^2-1)/2,v^4(v^2-2)/4)$-cube $C_1$.
By \cite{BArxiv22}, there exists a layer-rainbow cube~$L_2$ of order $v^2$
containing~$L_1$, with entries in $\{1,\ldots,v^4\}$. From  Theorem~\ref{tmhadsqconstr}~(i),
when applied to~$K$ and~$L_2$, a $\C^3_3(v^6,v^8(v^4-1)/2,v^8(v^4-2)/4)$-cube $C_2$
is obtained. Choosing the embedded coordinates of~$H$ in~$K$ and the copy $L_1\subseteq L_2$ compatibly, the restriction of~$C_2$ to those coordinates is exactly~$C_1$.
\end{proof}

\subsection{Constructions using association schemes on triples}

Association schemes on triples were
introduced by Mesner and Bhattacharya~\cite{MB90, MB94}.
Let $\Omega$ be a non-empty finite set. Let $m\geq 4$, and let $\{R_0,R_1,\ldots,R_m\}$ be
a partition of $\Omega\times \Omega\times \Omega$. For a permutation $\sigma$ of $\{1,2,3\}$, let
$\sigma(R_i)=\{(x_{\sigma(1)},x_{\sigma(2)},x_{\sigma(3)}) \mid (x_1,x_2,x_3)\in R_i\}$. The pair
$\big(\Omega,\{R_i\}_{i=0}^m\big)$ is an \emph{association scheme on triples} --- abbreviated as \emph{AST} --- if the following  conditions hold.
\begin{enumerate}[label=(\Roman*)]
    \item\label{ast1} The relations $R_0,R_1,R_2,R_3$ are chosen so that
 \begin{align*}
    &&
    &R_0=\{(x,x,x)\mid x\in \Omega\},
    &&
    R_1=\{(x,y,y)\mid x,y\in \Omega,x\neq y\},\\
    &&&
    R_2=\{(x,y,x)\mid x,y\in \Omega,x\neq y\},
    &&
    R_3=\{(x,x,y)\mid x,y\in \Omega,x\neq y\}.
    &&
  \end{align*}

    \item\label{ast2} For any $i\in\{0,1,\ldots,m\}$ and any permutation $\sigma$ of $\{1,2,3\}$,
    $\sigma(R_i)=R_j$ holds for some $j\in\{0,1,\ldots,m\}$.
    \item\label{ast3} For any two distinct elements $y,z\in \Omega$ and $i\in \{0,1,\ldots,m\}$, the size of the set $\{x\in \Omega \mid (x,y,z)\in R_i\}$ depends only on $i$, and not on the choice of $y$, $z$.
    \item\label{ast4} For any $i,j,k,\ell\in\{0,1,\ldots,m\}$ and any triple $(x,y,z)\in R_\ell$, the size of the set
    \[
    \{w\in \Omega \mid (w,y,z)\in R_i,(x,w,z)\in R_j,(x,y,w)\in R_k\}
    \]
    depends only on $i,j,k,\ell$, and not on the choice of $(x,y,z)$.
    This value is said to be the \emph{intersection number} and denoted by $p_{ijk}^{\ell}$.
\end{enumerate}
The trivial relation $R_0$ from~\ref{ast1} is the support of the identity cube from
Proposition~\ref{propidcube}. Regarding condition~\ref{ast3}, we write
\[
n_i^{(1)}=|\{x\in \Omega \mid (x,y,z)\in R_i\}|
\]
for distinct $y,z\in \Omega$. Similarly, we define
\[
n_i^{(2)}=|\{y\in \Omega \mid (x,y,z)\in R_i\}|, \quad
n_i^{(3)}=|\{z\in \Omega \mid (x,y,z)\in R_i\}|.
\]
Because of condition~\ref{ast2}, these numbers also do not depend on the choice of
distinct elements $x,z\in \Omega$ and $x,y\in \Omega$, respectively. It is obvious
that $n_0^{(i)}=n_i^{(i)}=0$ and $n_j^{(i)}=1$ for $j\in\{1,2,3\}\setminus \{i\}$.

We say that $R_i$ is \emph{symmetric} with respect to the permutation $\sigma$ if $\sigma(R_i)=R_i$.
We call the AST \emph{symmetric} if each relation $R_i$,
$i\in\{4,\ldots,m\}$, is symmetric for any permutation~$\sigma$ of $\{1,2,3\}$. In that case,
$n_i^{(1)}=n_i^{(2)}=n_i^{(3)}$ holds and we denote the common value by~$n_i$. For convenience, set $n_i:=n_i^{(3)}$ also for $i\in\{0,1,2,3\}$.

If $|\Omega|=v$, the AST is said to be of \emph{order}~$v$. To each relation
$R_i$, we associate a $v\times v\times v$ matrix~$A_i$ with entries
\[
A_i(x,y,z)=\begin{cases}
    1 & \text{ if }(x,y,z)\in R_i,\\
    0 & \text{ otherwise}.
\end{cases}
\]
The matrix $A_i$ is the \emph{adjacency matrix} of the hypergraph $(\Omega,R_i)$.
Conversely, the relation~$R_i$ is the support of the matrix~$A_i$.
Given three $v\times v\times v$ matrices $A$, $B$, and $C$, the \emph{ternary product}
$D=ABC$ is the $v\times v\times v$ matrix with entries
\[
D(x,y,z)=\sum_{w\in \Omega}A(w,y,z)B(x,w,z)C(x,y,w).
\]
Using this product, condition~\ref{ast4} in the definition of ASTs can
be expressed as
$$A_i A_j A_k =\sum_{\ell=0}^m p_{ijk}^\ell A_\ell, \quad \mbox{ for  }i,j,k\in \{0,\ldots,m\}.$$

ASTs are used in~\cite{BS} for a construction of Hadamard
cubes of propriety~$3$. The following theorem gives a similar
construction of $\C_3^3(v,k,\lambda)$-cubes.

\begin{theorem}\label{thm:ast}
Let $(\Omega,\{R_i\}_{i=0}^m)$ be a symmetric AST of order~$v$, and let $I\subseteq\{0,1,\ldots,m\}$ be a nonempty subset with either
$I\supseteq\{1,2,3\}$, or $I\cap \{1,2,3\} = \emptyset$. Then, the
$v\times v\times v$ array
$A=\sum_{i\in I} A_i$ is a $\C^3_3(v,k,\lambda)$-cube, where
\[
k:=\sum_{i,j\in I}\sum_{r=0}^m (p_{ijr}^0+(v-1)p_{ijr}^3), \quad
\lambda:=\sum_{i,j\in I}\sum_{r,\ell=0}^m  p_{ijr}^\ell n_\ell.
\]
\end{theorem}

\begin{proof}
By the condition on the subset~$I$ and because the AST is symmetric,
$A(x,y,z)=A(x',y',z')$ holds for any permutation $(x',y',z')$
of $(x,y,z)$. Therefore, it is enough to compute
$$\sum_{z,w\in \Omega}A(x,z,w)A(y,z,w).$$ Indeed, we calculate $AAJ$ in two ways.
On the one hand, using the property that $A(x,y,z)=A(y,x,z)$,
\begin{align}
        (AAJ)(x,y,z)&=\sum_{w\in \Omega}A(w,y,z)A(x,w,z)J(x,y,w) \notag\\
        &=\sum_{w\in \Omega}A(w,y,z)A(x,w,z)
        =\sum_{w\in \Omega}A(y,w,z)A(x,w,z).\label{eq:1}
    \end{align}
On the other hand,
    \begin{align*}
        AAJ=\sum_{i,j\in I}\sum_{r=0}^m A_i A_jA_r=\sum_{i,j\in I}\sum_{r,\ell=0}^m  p_{ijr}^\ell A_\ell,
    \end{align*}
    which shows
    \begin{align}
        (AAJ)(x,y,z)=\sum_{i,j\in I}\sum_{r,\ell=0}^m p_{ijr}^\ell A_\ell(x,y,z).\label{eq:2}
    \end{align}
Combining \eqref{eq:1} and \eqref{eq:2} and taking the sum over $z\in \Omega$,
\begin{align*}
\sum_{z,w\in \Omega}A(y,w,z)A(x,w,z)&=\sum_{i,j\in I}\sum_{r,\ell=0}^m p_{ijr}^\ell \sum_{z\in \Omega}A_\ell(x,y,z)
=\begin{cases}\displaystyle\sum_{i,j\in I}\sum_{r,\ell=0}^m p_{ijr}^\ell n'_\ell & \text{ if }x=y,\\[2mm]
\displaystyle\sum_{i,j\in I}\sum_{r,\ell=0}^m p_{ijr}^\ell n_\ell & \text{ if }x\neq y,
\end{cases}
\end{align*}
where $n'_{\ell}=|\{z\in \Omega \mid (x,x,z)\in R_\ell\}|=\begin{cases}1 & \text{ if }\ell=0,\\ v-1 & \text{ if }\ell=3, \\ 0 & \text{ otherwise}.\end{cases}$
\end{proof}

By Theorem~\ref{thm:ast}, the union of trivial relations $R_1\cup R_2\cup R_3$,
as well as each of the remaining relations of the AST, supports a $\C_3^3$-cube
indexed by~$\Omega$. Together, all these relations partition the Cartesian
product $\Omega\times \Omega\times \Omega$. Analogous two-dimensional objects
were named \emph{mosaics of combinatorial designs} in~\cite{GGP18}.
Additionally, the union of any subset of these relations also
supports a $\C_3^3$-cube. Mosaics with this remarkable property
were studied in~\cite{MK17, MSMKK06, SMMKK07} under the name
\emph{additive BIBDs} (not to be confused with the additive block
designs of~\cite{CFP17}).

\begin{corollary}\label{trivrel}
For any $v\ge 2$, the $v\times v\times v$ array $A$ defined by
$$A(x,y,z)=[x=y\neq z \mbox{ or } x=z\neq y \mbox{ or } y=z\neq x]$$
is a $\C^3_3(v,3(v-1),v)$-cube.
\end{corollary}

\begin{proof}
Consider the trivial association scheme on triples $(\Omega,\{R_i\}_{i=0}^4)$ with $|\Omega|=v$ and
take $A=\sum_{i=1}^3 A_i$ in Theorem~\ref{thm:ast}.
\end{proof}

In~\cite{MB90}, non-trivial families of ASTs are constructed
from Steiner $2$-designs, $2$-transitive groups, and regular
two-graphs. Some new constructions are obtained
in~\cite{BS, BB25, BP21}, and a recent survey is~\cite{BB23}.

\begin{remark}\textup{
By~\cite[Theorem~3.1]{MB90}, a Steiner system $S(2,k,v)$ gives rise
to a symmetric AST of order~$v$ with $m=5$. For distinct points~$x$, $y$, and~$z$,
the triple $(x,y,z)$ is in the relation~$R_4$ if the points are collinear,
and in~$R_5$ otherwise. Spence~\cite{ES96} proved
that there are exactly~$18$ Steiner systems $S(2,4,25)$ up to isomorphism.
The corresponding adjacency matrices~$A_4$ are $18$ nonequivalent
$\C^3_3(25, 48 , 2)$-cubes, and the matrices~$A_5$ are $18$ nonequivalent
$\C^3_3(25,504,420)$-cubes. Indeed, for either relation a line obtained by fixing two equal coordinates is zero, whereas fixing two distinct coordinates gives a positive line sum; hence any equivalence forces the three coordinate permutations to coincide, and therefore induces an isomorphism of the Steiner systems.
}\end{remark}

\section{Concluding remarks and research directions}\label{sec5}

\begin{figure}[!b]
\ifthenelse{\boolean{images}}
{\ifthenelse{\boolean{lowres}}
{\includegraphics[width=6cm]{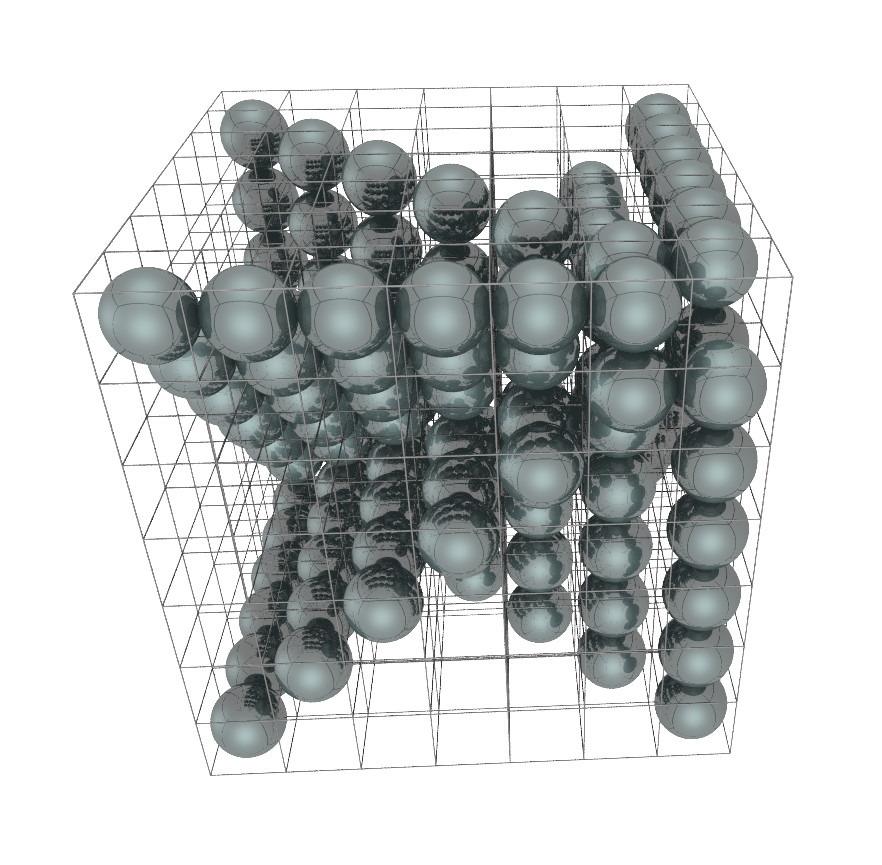}}
{\includegraphics[width=6cm]{7-18-7.jpg}}
}{}
\vskip -3mm
\caption{The cube of Corollary~\ref{trivrel} for $v=7$.}\label{figx}
\end{figure}

Consider the $\C^3_3(v,3(v-1),v)$-cube~$A$ of Corollary~\ref{trivrel}.
By the definition, any layer of~$A$ contains exactly $k=3(v-1)$ entries~$1$,
and any two parallel layers contain exactly $\lambda=v$ common \hbox{$1$-entries}.
In addition, this cube has the following property: Any~$t$ parallel layers contain exactly
$\lambda_t=v-t$ common $1$-entries, for $t\in \{3,\ldots,v\}$.
If we include the relation $R_0$, i.e.\ add the
diagonal to~$A$, a $\C^3_3(v,3v-2,v+2)$-cube with the same
property for $\lambda_t=v$, $t\in \{3,\ldots,v\}$ is obtained.
The cube~$A$ is depicted in Figure~\ref{figx} for $v=7$.

This example is reminiscent of combinatorial $t$-$(v,k,\lambda)$ designs,
where every set of~$s$ points is contained in exactly~$\lambda\binom{v-s}{t-s}/\binom{k-s}{t-s}$ blocks, for any $s\leq t$~\cite{KL07}. The largest such integer~$t$ is
called the \emph{strength} of the design. It is well known that the strength
of nontrivial $2$-dimensional symmetric designs is bounded by $t\le 2$,
but in view of the example it makes sense to study $n$-dimensional symmetric
designs of higher strength. We propose the following direction for further
research.

\begin{problem}
Study symmetric $(v,k,\lambda)$ block designs of dimension~$n$,
propriety~$d$, and strength $t$.
\end{problem}

Two-dimensional symmetric $(v',k',\lambda')$ designs with
$(v,k,\lambda)$ subdesigns have been extensively studied; see
the survey~\cite{MS02}. \'{O}~Cath\'{a}in~\cite{POC13} proved that
the case $v'=v$, $k'=k+1$ is equivalent to doubly regular tournaments.
Theorem~\ref{embsymcub} provides a family of $\C^3_3(v',k',\lambda')$-cubes
with subcubes in $\C^3_3(v,k,\lambda)$. The cubes are of order $v'=h^6$, and the
subcubes of order $v=h^3$, where $h^2$ is the order of a Hadamard
matrix. Examples of cubes with subcubes of the same order $v'=v$ are obtained
from Theorem~\ref{thm:ast}: take any two sets $I\subseteq I'\subseteq \{0,\ldots,m\}$
satisfying the condition of the theorem. If $0\notin I$ and we take $I'=I\cup\{0\}$, then
$k'=k+1$ also holds. We propose to study this problem further for three-dimensional
symmetric designs.

\begin{problem} \label{embsym3d}
Find conditions under which $\C^3_3(v',k',\lambda')$-cubes with
$\C^3_3(v,k,\lambda)$-subcubes exist, for $v\le v'$, $k\le k'$,
and $\lambda\le \lambda'$.
\end{problem}

In particular, in relation to Theorem~\ref{embsymcub}, we put forward the following conjecture.
\begin{conjecture}
Let $H$ be a regular Hadamard matrix of order $v^2$. Then there exists a $\C^3_3(v^6,v^6(v^6-1)/2,v^6(v^6-2)/4)$-cube containing a $\C^3_3(v^3,v^3(v^3-1)/2,v^3(v^3-2)/4)$-subcube.
\end{conjecture}
The existence follows from Theorem~\ref{tmhadsqconstr}~(ii)
if~$H$ can be embedded in a regular Hadamard matrix of order~$v^4$, but
we are not aware of any results on embedding regular Hadamard matrices.
In connection with Problem \ref{embsym3d} and Ryser's embedding question for
Hadamard matrices~\cite{Ne82, Wal86}, we propose the following problem.
\begin{problem}\label{proemb3dabc}
Given an $a\times b\times c$ array $A$ over $\{0,1\}$, determine the smallest possible order $f(A)$ of a three-dimensional symmetric design containing $A$.
\end{problem}
It would be interesting if one could produce such results to obtain other families of
three-dimensional symmetric designs containing subdesigns.

In this paper, we have mainly considered the existence problem for
$\C^3_3(v,k,\lambda)$-cubes. It appears that the enumeration problem
is exceedingly difficult, but in some cases it may be possible to
get upper or lower bounds on the number of cubes.

\begin{problem}
     Find good bounds for $|\C^{n}_d(v,k,\lambda)|$.
\end{problem}

For example, for $n\ge 2$, Linial and  Luria~\cite{LL14} defined an $n$-dimensional
permutation of order~$v$ as a $v\times \dots \times v$ array of
dimension~$n+1$ over $\{0,1\}$ in which every row contains a unique
$1$-entry. These objects are precisely $\C^{n+1}_3(v,v,0)$-cubes,
and they are equivalent to $n$-dimensional permutation
cubes of~\cite{KD15}, or $n$-dimensional Latin hypercubes
of~\cite{MW08}. The equivalence is a higher-dimensional
variant of Proposition~\ref{propls}. In~\cite{LL14},
it is conjectured that $|\C^{n+1}_3(v,v,0)|=
\big((1+ o(1))\dfrac{v}{e^n}\big)^{v^n}$ and it is proved that the upper
bound holds.

%\textcolor{red}{I think we can count $|\C^{3}_3(v,k,0)|$ for $k\leq v$.}

%\begin{comment}
\subsubsection*{Acknowledgments}
Amin Bahmanian's research is partially supported by a  Faculty Research Award at Illinois State University. Sho Suda's research is supported by JSPS KAKENHI Grant Number 22K03410. Vedran Kr\v{c}adinac and Lucija Reli\'{c} are partially supported by the European Union -- NextGenerationEU through the National Recovery and Resilience Plan 2021-2026 Institutional grants of University of Zagreb Faculty of Science (IK IA 1.1.3 Impact4Math).
%\end{comment}

\end{document}